\documentclass[leqno,12pt]{article}
\usepackage{amsmath,amssymb}

\newtheorem{Theorem}{\hspace{\parindent}\bf Theorem}[section]
\newtheorem{Lemma}{\hspace{\parindent}\bf Lemma}[section]
\newtheorem{Proposition}{\hspace{\parindent}\bf Proposition}[section]

\newtheorem{Definition}{\hspace{\parindent}\bf Definition}[section]

\newcommand{\qed}{\hfill$\square$\vspace{0.3cm}}
\newcommand{\Tr}{\mathop{\rm Tr}}
\newcommand{\E}{\mathop{\rm E}}

\begin{document}

\title{\textbf{A fractional porous medium equation }}
\author{ by \\
Arturo de Pablo, Fernando Quiros, \\ Ana Rodriguez and Juan Luis Vazquez}

\maketitle

\begin{abstract}
We develop a theory of existence, uniqueness and regularity for the following
porous medium equation  with fractional diffusion,
$$
\left\{
\begin{array}{ll}
\dfrac{\partial u}{\partial t} + (-\Delta)^{1/2} (|u|^{m-1}u)=0, &
\qquad  x\in\mathbb{R}^N,\; t>0,
\\ [8pt]
u(x,0) = f(x), & \qquad x\in\mathbb{R}^N,%
\end{array}
\right.
$$
with $m>m_*=(N-1 )/N$, $N\ge1$ and $f\in L^1(\mathbb{R}^N)$. An
$L^1$-contraction  semigroup is constructed and the continuous dependence on
data and exponent is established. Nonnegative solutions are proved to be
continuous and strictly positive for all $x\in\mathbb{R}^N$, $t>0$.

\end{abstract}

\section{Introduction}
\setcounter{equation}{0}

This paper is concerned with the existence, uniqueness  and properties of
solutions $u=u(x,t)$ to the Cauchy problem
\begin{equation}  \label{eq:main}
\left\{
\begin{array}{ll}
\dfrac{\partial u}{\partial t} + (-\Delta)^{1/2} (|u|^{m-1}u)=0, & \qquad  x\in\mathbb{R}^N,\; t>0,
\\ [4mm]
u(x,0) = f(x), & \qquad x\in\mathbb{R}^N,%
\end{array}
\right.
\end{equation}
for  exponents $m>0$, in space dimension $N\ge 1$,  and with initial value
$f\in L^1(\mathbb{R}^N)$.  By a solution it is meant a suitable concept of weak
or strong solution. In particular, we prove that $u \in C( [0,\infty): L^1
(\mathbb{R}^N))$ and that the equation is satisfied a.e.~in
$Q=\mathbb{R}^N\times(0,\infty)$. The sign requirement $u\ge 0$ is not strictly needed but when enforced some additional  properties hold.

We recall that the nonlocal operator
$(-\Delta)^{1/2}$ is defined  for any function $g$ in the  Schwartz class
through the Fourier transform,
\begin{equation}
\widehat{(-\Delta)^{1/2}  g}\,(\xi)=|\xi| \,\widehat{g}(\xi),
\label{def-fourier}
\end{equation} or via the Riesz potential,
\begin{equation}
(-\Delta)^{1/2}  g(x)= C_{N }\mbox{ P.V.}\int_{\mathbb{R}^N}
\frac{g(x)-g(y)}{|x-y|^{N+{1} }}\,dy,
\label{def-riesz}
\end{equation}
where $C_{N }=\pi^{-\frac{N+1}2}\Gamma(\frac{N+1}2)$ is a normalization
constant, see for example \cite{Landkof}, \cite{Stein}.

Equations of this form can be considered as nonlinear variations of the linear
fractional diffusion equation obtained for $m=1$, which is a model of so-called
anomalous diffusion, a much studied topic in physics, probability and finance,
see for instance \cite{AT, Jara2, JKOlla, MMM, VIKH, WZ} and the references
therein. We recall that fractional Laplacian operators of the form
$(-\Delta)^{\sigma/2}$, $\sigma\in(0,2)$, are infinitesimal generators of
stable L{\'e}vy processes
\cite{Applebaum,Bertoin}. The analysis of the linear equation in the
whole space is easy since an integral representation can be used for
the solutions, see below. Such a representation is not available in
the nonlinear case.

Interest in studying the nonlinear model we propose is two-fold: on the one
hand, experts in the mathematics of diffusion want to understand the
combination of fractional operators with porous medium type propagation, and on
the other hand models of this kind  arise in statistical mechanics
\cite{Jara} and heat control \cite{AC}. The rigorous study of such nonlinear models has
been delayed by mathematical difficulties in treating at the same time the
nonlinearity and fractional diffusion.

Observe that the above equation becomes the well-known Porous Medium Equation
(PME) when replacing the  nonlocal diffusion operator $(-\Delta)^{1/2}$ by the classical
Laplacian $-\Delta$. A number of techniques in dealing with the present
nonlinear fractional diffusion model will be borrowed from the experience
obtained with the PME, see for instance \cite{Vazquez:Book}. Our original
purpose was to study Problem \eqref{eq:main} for every $m> 1$, to examine the
existence and properties of \lq\lq fractional slow diffusion". But the
development of the theory allows to cover with a reasonable additional effort
the ``fast diffusion cases", $m<1$, on the condition that we restrict the
exponent to be larger than a critical value, $m>m_*\equiv(N-1 )/N$. This
critical value is intrinsic to the equation, it appears in various contexts of
the theory. It corresponds to the classical critical value $m_*=(N-2)_+/N$ in
the PME case, see
\cite{Benilan-Crandall}. Existence of a weak solution is however proved for
every $m>0$.

\noindent {\sc Harmonic extensions. } Besides formulae \eqref{def-fourier} and
\eqref{def-riesz}, there is another way of computing the half Laplacian, through the so-called Dirichlet to
Neumann operator. If $g=g(x)$ is a smooth bounded function defined in
$\mathbb{R}^N$, we consider its harmonic extension $v=v(x,y)$ to the upper
half-space $\mathbb{R}^{N+1}_+$, $v=\E(g)$, i.\,e.,  the unique smooth bounded
solution to
\begin{equation}
\left\{
\begin{array}{ll}
\Delta_{x,y}  v=0,\qquad &x\in\mathbb{R}^N,\, y>0,\\
v(x,0)=g(x),\qquad&x\in\mathbb{R}^N.
\end{array}
\right. \label{alfa-extension}\end{equation}
Then,
\begin{equation}
-\frac{\partial v}{\partial y}(x,0)=(-\Delta_x)^{1/2}  g(x)\,,
\label{fract-lapla}
\end{equation}
where  $\Delta_{x,y}$ is the Laplacian in all $(x,y)$ variables and $\Delta_x$
acts only on the $x$ variables (in the sequel we will drop the subscripts when
no confusion arises). In order to check~\eqref{fract-lapla}, just apply the
operator in the right-hand side twice.  The choice of sign for the normal
derivative makes the operator positive. Observe that the extension operator is
well defined in $H^{1/2}(\mathbb{R}^N)$, and so is the Dirichlet to Neumann
operator, which coincides with $(-\Delta)^{1/2}$ in this more general setting.
This well known technique has been recently used in several situations, see for
instance
\cite{Cabre-Tan,CS,Silvestre}.

\noindent {\sc Problem-setting. } By means of the above-mentioned harmonic extension
we rewrite, for smooth solutions, the nonlocal Problem \eqref{eq:main} in a ``local
way'' (i.\,e., using local differential operators) as a quasi-stationary problem with a
dynamical boundary condition.  Indeed, $w=|u|^{m-1}u$ satisfies
\begin{equation}
\left\{
\begin{array}{ll}
\Delta   w=0\qquad &\mbox{for } x\in\mathbb{R}^N,\,y>0,\, t>0,\\
\dfrac{\partial w}{\partial y}-\dfrac{\partial (|w|^{\frac1m-1}w)}{\partial
t}=0\qquad&\mbox{for } x\in\mathbb{R}^N,\,y=0,\, t>0,\\
w(x,0,0)=f^m(x)\qquad&\mbox{for } x\in\mathbb{R}^N.
\end{array}
\right. \label{pp:local0}
\end{equation}
This problem has  been recently considered  by Athanasopoulos and Caffarelli
\cite{AC}. They prove that  bounded weak energy solutions to \eqref{pp:local0}
are H{\"o}lder continuous if $m>1$. The existence and uniqueness of that kind of
solutions is one of the outcomes of the present paper.

The connection between problems with dynamical boundary conditions
and nonlocal equations has already been exploited in \cite{Viti2} in
the case of a bounded domain, and in \cite{AF} for a semilinear
problem in the half-space $\mathbb{R}^{N+1}_+$. However, in those
works the study of the nonlocal equation is used to obtain
properties of the local one. Here, our approach is exactly the
opposite.

\medskip

\noindent {\sc Main results. }  Our purpose is to establish a theory of existence, uniqueness,
comparison and regularity for suitable weak solutions of Problem
\eqref{eq:main} with initial data $f\in L^1(\mathbb{R}^N)$. The full theory
works for values of $m$ larger than the critical value $m_*$ mentioned above,
but basic existence and uniqueness holds for all $m>0$, for data which are
moreover bounded.

 Section \ref{sect-prelim} contains preliminaries,
the basic definitions of solutions, and a list of main results.
We define the concept of weak solution to Problem \eqref{eq:main}
through the standard concept of weak solution to the
associated local Problem \eqref{pp:local0}. We also define the concept of strong solution.

We establish the existence of weak solutions in Section \ref{sect-existence}
by means of Semigroup Theory, solving first some associated
elliptic problem, under the
condition that the initial data $f$ are both integrable and bounded. Actually, the obtained solution is strong  and the equation is satisfied almost
everywhere. We also prove in this context an $L^1$-$L^\infty$ estimate that will be basic for the so-called smoothing effect.

Uniqueness is studied in Section~\ref{sect-uniq}. Section~\ref{sect-props} deals with further  properties of the constructed solution. It includes conservation of mass, positivity and regularity. At this point we use the continuity result from \cite{AC} to show that solutions to \eqref{eq:main} corresponding to nonnegative initial data become immediately
strictly positive if $m>1$. This is a remarkable property since it  departs
from the well-known properties of the standard  PME, cf. \cite{Vazquez:Book}.
On the other hand, for $m_*<m<1$ we are able to prove the expected positivity
property using a different approach. This is used later, in combination with
boundedness and a result in \cite{AC} to prove H{\"o}lder continuity also in this
case. Let us notice that, unlike in the local case, there is still no general
regularity result for linear nonlocal equations (with reasonable coefficients)
guaranteeing that positive bounded solutions to problem \eqref{eq:main} are in
fact $C^\infty$, though this is expected to be true.

After such a work, we are able  to treat
general solutions with data in $L^1(\mathbb{R}^N)$ in Section~\ref{sect-L1}.
Here we complete the proof of uniform boundedness of the solutions with integrable data for
positive times, the  $L^1$-$L^\infty$ smoothing effect.

In Section~\ref{sect.cont.dep.} we study the continuous dependence of the solution in terms of the
exponent $m$ and the data $f$, in the case $m>m_*$. In particular, we show that the linear case $m=1$ can be obtained as a limit of the nonlinear case both from above and below.

Section \ref{sec.cae} contains a brief description of alternative approaches to the existence theory and an announcement of extensions. Finally, the Appendix gathers some  technical lemmas.
\medskip

\noindent {\sc  Notice on the linear case.} For the value of the parameter $m=1$ we obtain
the equation
\begin{equation}
\frac{\partial u}{\partial t}+(-\Delta)^{1/2}u=0.
\label{linearcase}\end{equation}
This is a linear fractional heat equation where the fractional derivatives act
only on the space variable. It is explicitly solvable in terms of the initial
value, $u(x,0)=f(x)$, through convolution with the explicit Poisson kernel in
$\mathbb{R}^{N+1}_+$,
\begin{equation}
u(x,t)=C_N\int_{\mathbb{R}^N}\frac{t\,f(z)}{(|x-z|^2+t^2)^{(N+1)/2}}\,dz,
\label{linear-rep}\end{equation}
where $C_N$ is the constant in \eqref{def-riesz}. Note that this corresponds to
an anomalous diffusion law of the form $ \langle x\rangle\sim t^{\alpha}$ with
$\alpha=1$ instead of the standard $\alpha=1/2$ of the Brownian case.

\medskip

\noindent {\sc Notations.} In dealing with extended functions, we
denote the upper half-space, $\mathbb{R}^{N+1}_+$, by $\Omega$, and write its
points as $\overline x=(x,y)$, $x\in\mathbb{R}^N$, $y>0$. We denote by $\Gamma$
the boundary of $\Omega$, i.\,e., $\Gamma=\mathbb{R}^N\times\{0\}$, which is
identified to the original $\mathbb{R}^N$  with variable $x$. We consider also
the extension and trace operators, $\E$, $\Tr$: for a function $v\in
H^{1/2}(\Gamma)$, we denote its harmonic extension to $\Omega$ as $\E(v)$;
notice that $\E(v)\in H^1(\Omega)$; on the other hand, given a function $z\in
H^1(\Omega)$, we denote its trace on $\Gamma$, which belongs to
$H^{1/2}(\Gamma)$, as $\Tr(z)$.

As in the  PME theory, we will be mostly interested in nonnegative data and
solutions. However, the basic theory can be developed for data of any sign, and
in that case we will use the simplified notation $u^m$ instead of the ``odd
power'' $|u|^{m-1}u$, and we will also use such a notation when $m$ is replaced
by $1/m$. In fact, we will show that if the initial value is nonnegative, then
the weak solution we construct is also nonnegative, $u\ge0$, which helps
justify our abbreviated notations.

\section{Preliminaries and main results}\label{sect-prelim}
\setcounter{equation}{0}

As mentioned above, we define the concept of weak solution to Problem
\eqref{eq:main} through the standard concept of weak solution to an associated
local problem, which we write here again by convenience.
\begin{equation}
\left\{
\begin{array}{ll}
\Delta   w=0\qquad &\mbox{for } \overline x\in\Omega,\, t>0,\\
\dfrac{\partial w}{\partial y}-\dfrac{\partial w^{1/m}}{\partial
t}=0\qquad&\mbox{on } \Gamma,\,
t>0,\\
w(x,0,0)=f^m(x)\qquad&\mbox{on } \Gamma.
\end{array}
\right. \label{pp:local}
\end{equation}
In order to define a weak solution of this problem we multiply formally the
equation in
\eqref{pp:local} by a test function $\varphi$ and integrate by parts to obtain
\begin{equation}
\displaystyle -\int_0^T\int_{\Omega}\langle\nabla w,\nabla
\varphi\rangle\,d\overline xds+\int_0^T\int_{\Gamma}u\dfrac{\partial
\varphi}{\partial t}\,dxds=0\,,
\label{weak-local}\end{equation}
with $u=\Tr(w^{1/m})$, on the condition that $\varphi$ vanishes for $t=0$ and
$t=T$, and also for large $|x|$ and $y$.

\begin{Definition}
We say that a pair of functions $(u,w)$ is a {\it weak solution} to Problem
\eqref{pp:local} if $w\in L^1((0,T);W^{1,1}_{loc}(\Omega))$,
$u=\Tr(w^{1/m})\in L^{1}(\Gamma\times (0,T))$ and equality
\eqref{weak-local} holds for every $\varphi\in C_0^1(\overline\Omega\times[0,T))$.
Finally, the initial data are taken in the sense that \ $\lim\limits_{t\to0}u(\cdot,t)=f$ in $L^1(\mathbb{R}^N)$.
\end{Definition}

An alternative form of equality \eqref{weak-local}, including the initial value
in it, is
\begin{equation}
\begin{array}{l}
\displaystyle-\int_0^T\int_{\Omega}\langle\nabla w,\nabla
\varphi\rangle\,d\overline xds+\int_0^T\int_{\Gamma}u\dfrac{\partial
\varphi}{\partial t}\,dxds=\\ [4mm]
\displaystyle\int_{\Gamma}u(x,T)\varphi(x,0,T)\,dx-
\int_{\Gamma}f(x)\varphi(x,0,0)\,dx.
\end{array}
\label{weak-local2}\end{equation}

As is usual,  more general test functions can be considered by approximation,
whenever the integrals make sense. Note that the trace $u=\Tr(w^{1/m})$ is  well
defined. For brevity we will refer sometimes to the solution as only $u$, or
even only $w$, when no confusion arises, since it is clear how to complete the
pair from one of the components, $u=\Tr(w^{1/m})$, $w=\E(u^m)$. By a weak solution
of our original Problem
\eqref{eq:main} we understand $u$, the first element of the solution to Problem~\eqref{pp:local}.

Observe that the definite advantage of working with the local version is
compensated in some sense by the difficulty of having integrals in
\eqref{weak-local} defined in spaces of different dimensions.

This definition is a very general notion of solution: in this framework we can
construct a weak solution to Problem \eqref{pp:local} provided the initial
value $f$ is integrable and bounded. We restrict ourselves in the next results
to such data. However, weak solutions are sometimes difficult to
work with, and we are not able to prove uniqueness. Hence, a class of solutions
with better properties is welcome. A  quite convenient choice is the class of
so-called {\it weak energy solutions}, cf. \cite{Vazquez:Book} for the standard
 PME.

\begin{Definition}
  A weak solution pair $(u,w)$ to Problem \eqref{eq:main} is said to be a weak energy solution
  if moreover $w\in L^2([0,T]; H^1(\Omega))$.
\end{Definition}

\begin{Theorem}\label{th:main} Let $m>0$. For every $f\in L^1(\mathbb{R}^{N})\cap L^
\infty(\mathbb{R}^{N})$ there exists a unique weak energy solution  to Problem
\eqref{pp:local}. Moreover $u\in C([0,\infty):\, L^1(\mathbb{R}^{N}))\cap
L^\infty(\mathbb{R}^{N}\times[0,\infty))$.
\end{Theorem}

The importance of this class of solutions, besides having uniqueness, is that,
if we restrict to  nonnegative data and exponents $m>m_*=(N-1)/N$, we can
obtain regularity and positivity.

\begin{Theorem}\label{th-prop-energy}
Let $f\in L^1(\mathbb{R}^N)\cap L^\infty(\mathbb{R}^N)$ be nonnegative, and
assume $m>m_*$. Then the weak energy solution $(u,w)$ to Problem
  \eqref{pp:local} satisfies:

  \noindent $i)$ Conservation of mass: for every $t>0$ we have
\begin{equation}
\int_{\mathbb{R}^N}u(x,t)\,dx=\int_{\mathbb{R}^N}f(x)\,dx\,.
\label{conservation}\end{equation}

\noindent $ii)$ Positivity:  $u(\cdot,t)>0$ in
$\mathbb{R}^N$for every $t>0$.

\noindent $iii)$ Regularity: there exists some $0<\alpha<1$ such
that $u\in C^\alpha(\mathbb{R}^N\times(0,T))$.

\noindent (iv) Maximum Principle: if $u_1, u_2$ are solutions with data $u_{01}, u_{02}$
and $u_{01}\le u_{02}$ a.e. in $\mathbb{R}^N$, then $u_{1}\le u_{2}$ a.\,e., in \
$Q=\mathbb{R}^N\times (0,\infty)$.

\noindent (v) Contraction: for any two solutions $u_1, u_2$ with
data $u_{01}, u_{02}$ we have
\begin{equation}
\|u_1(\cdot,t)-u_2(\cdot,t)\|_{L^1(\mathbb{R}^N)}\le \|u_{01}-u_{02}\|_{L^1(\mathbb{R}^N)}.
\end{equation}
\end{Theorem}

The restriction $m>m_*$ is not technical: positivity and conservation of mass
are not true if $m<m_*$, see Proposition~\ref{pro:extinction}. On the other
hand, conservation of mass holds also for solutions with changing sign if
$m>m_*$. If $m>1$ the $C^\alpha$ regularity result is true also for any
changing sign solution, \cite{AC}.

A further interesting property is that the weak energy solutions are
{\sl strong solutions}, which means that the terms (in principle
only distributions) involved in equation \eqref{pp:local} are in
fact functions, and equalities hold almost everywhere. The main
technical difficulty is to prove that $\partial_t u$ is a function.

\begin{Theorem} In the hypotheses of Theorem {\rm \ref{th-prop-energy}} we have $\partial_t u\in
L^1(\mathbb{R}^N)$.
\end{Theorem}

We observe that for strong solutions we can multiply \eqref{pp:local} by any
integrable function to get, instead of
\eqref{weak-local}, the following identity
\begin{equation}
\int_{\Omega}\langle\nabla w,\nabla
\varphi\rangle\,d\overline x +\int_{\Gamma}\frac{\partial
u}{\partial t}\,\varphi\,dx=0.
\label{strong-local}\end{equation}

Working with strong solutions we can use the solution itself as a test function
in formula \eqref{strong-local}. In particular this allows us to obtain a
universal bound for all solutions with the same mass.

\begin{Theorem} Let $f\in L^1(\mathbb{R}^N)\cap L^\infty(\mathbb{R}^N)$, and
assume $m>m_*$. Then,  there exists a positive constant $C$ such that the weak
energy solution to Problem  \eqref{pp:local} satisfies
\begin{equation}
\sup_{x\in\mathbb{R}^N}|u(x,t)|\le C\,t^{-\gamma }\|f\|_{L^1(\mathbb{R}^N)}^{\gamma/N}
\label{eq:L-inf}\end{equation}
with $\gamma=(m-1+{1} /N)^{-1}$. The constant $C$ depends only on $N$ and $m$.
\label{th:L-inf}\end{Theorem}

\noindent {\sc General integrable data.}  Once this theory is settled,  we are interested
in considering all integrable functions $f$ as possible data in Problem
\eqref{eq:main}. As we have advanced, this can be managed by approximation by
bounded initial data, and this is possible if we have an $L^1$-contraction at
hand. We thus introduce the concept of $L^1$ energy solution: a weak solution,
continuous in $L^1$, which is also an energy solution for positive times.

\begin{Definition} We say that a weak solution $(u,w)$  to Problem \eqref{pp:local} is an $L^1$-energy solution if
$u\in C([0,\infty):\, L^1(\mathbb{R}^{N}))$ and  $|\nabla w|\in
L^2(\Omega\times[\tau,\infty))$, for every $\tau>0$.
\end{Definition}

The $L^1$-contraction property for $L^1$-energy solutions is as follows.

\begin{Theorem}\label{th:contract} Let $(u,w)$ and $(\widetilde u,\widetilde w)$ be two $L^1$-energy solutions
to Problem \eqref{pp:local}. Then, for every $0\le t_1<t_2$,
\begin{equation}
\int_{\mathbb{R}^N}[u(x,t_2)-\widetilde u(x,t_2)]_+\,dx\le
\int_{\mathbb{R}^N}[u(x,t_1)-\widetilde u(x,t_1)]_+\,dx.
\label{eq:contract-l1}\end{equation}
\end{Theorem}

We therefore have that, when performing the approximation by problems with
bounded data, the limit function obtained is an $L^1$-energy solution. Now,
since estimate
\eqref{eq:L-inf} does not depend on the $L^\infty$ norm of the data, it is also true
for the limit solution (for changing sign solutions it holds by comparison). In
particular this represents an $L^1$-$L^\infty$ smoothing effect that allows to
obtain the same properties of Theorem \ref{th-prop-energy} for positive times.

\begin{Theorem}\label{th:main2} Let  $m>m_*$. Then for every $f\in L^1(\mathbb{R}^N)$ there exists
a unique $L^1$-energy solution to Problem \eqref{pp:local}. It satisfies estimate
\eqref{eq:L-inf} and the conservation of mass \eqref{conservation}. If moreover $f\ge0$, positivity and
regularity also hold, and the solution is strong. The maps $S_t:f\mapsto u(t)$
generate a nonlinear semigroup of order-preserving contractions in
$L^1(\mathbb{R}^N)$.
\end{Theorem}

The next sections  of the paper are devoted to treat the case of bounded data.
In Section~\ref{sect-L1} we drop this restriction and deal by approximation
with general $L^1$ data. We point out that the continuous dependence of the solutions constructed above with respect to the initial data and the exponent is stated and proved in Section \ref{sect.cont.dep.}.


\section{Weak solutions}\label{sect-existence}
\setcounter{equation}{0}

We set out to construct  a weak solution to the extended local Problem
\eqref{pp:local} taking initial values $f\in L^1(\Gamma)\cap L^\infty(\Gamma)$.
We point out that the construction of a weak solution can be made for every
$m>0$ and data not necessarily signed.

A well-known method of construction of solutions of evolution equations, and
also of generating a semigroup in a convenient functional space, is the
so-called Implicit Time Discretization. It runs as follows: if the evolution
equation is $dv/dt+A(v)=0$, where $A$ is a linear or nonlinear, bounded or
unbounded operator acting on a Banach space $\cal X$, and given initial data
$v(0)=f\in \cal X$, then the construction of an approximate solution of the
problem in a time interval $[0,T]$ proceeds dividing the time interval $[0,T]$
in $n$ subintervals of length $\varepsilon=T/n$ and then defining the
approximate solution $v_\varepsilon$ constant on each subinterval in the
following way: in each interval $(t_{k-1},t_k]$, $t_k=k\varepsilon$,
$k=1,\cdots,n$, we consider the solution $v_{\varepsilon,k}$ to the discretized
problem
\begin{equation}
\frac1{\varepsilon}(
v_{\varepsilon,k}-v_{\varepsilon,k-1})+A(v_{\varepsilon,k})=0.
\end{equation}
We take as starting condition $v_{\varepsilon,0}=f_{\varepsilon}$, an
approximation of $f$.  In the case of linear operators, a variant of the
Hille-Yosida Theorem ensures the convergence of these approximate solutions to
the so-called {\sl mild solution} of the evolution problem when the operator
$A$ satisfies some properties, like being maximal monotone, cf.~\cite{Brezis}.
The convergence result in the case of nonlinear and possibly unbounded
operators is given by the famous Crandall-Liggett Theorem~\cite{CL} under the
assumption that $A$ must be accretive and satisfy a certain rank condition.
(Reminder: a possibly nonlinear and unbounded operator $A:D(A)\subset {\cal
X}\to {\cal X}$ is called accretive if for every $\varepsilon>0$ the map
$I+\varepsilon A$ is one-to-one onto a subspace $R_\varepsilon (A) \subset
{\cal X}$ and the inverse ${\cal R}(\varepsilon, A):=(I+\varepsilon A)^{-1}:
R_\varepsilon (A)\to
\cal X$ is a contraction in the $\cal X$-norm. The precise rank condition that
we will use is $R_\varepsilon (A)\supset\overline{D(A)}$ for every $\varepsilon
>0$).

One of the typical examples of such theory is the standard  PME posed on the
whole space or on a bounded domain with homogeneous boundary conditions. The
early work due to B\'enilan and collaborators,
\cite{Benilan}, drew attention to this important results, as well as the
application to more general nonlinear diffusion-convection models.

We will apply such a strategy to our evolution Problem \eqref{pp:local}.  The
discretized problem is:

\begin{center}\begin{minipage}{14cm}Given $f\in L^1(\Gamma)\cap L^\infty(\Gamma)$ and
$\varepsilon>0$, to find $u_{\varepsilon}=\{u_{\varepsilon,1}, \cdots,
u_{\varepsilon,n}\}$ by solving for $k=1,\cdots, n$ the problem
\begin{equation}
\left\{
\begin{array}{ll}
\Delta w_{\varepsilon,k}=0\qquad &\mbox{ in } \Omega,\\
\varepsilon \dfrac{\partial w_{\varepsilon,k}}{\partial
y}=u_{\varepsilon,k}- u_{\varepsilon,k-1}\qquad&\mbox{ on }
\Gamma,
\end{array}
\right. \label{pe:local-iterate}
\end{equation}
with initial value $u_{\varepsilon,0}=f$ on $\Gamma$. In each such step,
$u_{\varepsilon,k-1}=\Tr(w_{\varepsilon,k-1}^{1/m})$ is known and
$u_{\varepsilon,k}$ and $w_{\varepsilon,k}=\E(u_{\varepsilon,k})$ are the
unknowns.
\end{minipage}\end{center}

The second equation in
\eqref{pe:local-iterate} can be written as
\begin{equation}
u_{\varepsilon,k}+\varepsilon A(u_{\varepsilon,k})=u_{\varepsilon,k-1},
\label{m-acc}\end{equation}
where the operator $A:D(A)\subset L^1(\Gamma)\to L^1(\Gamma)$ is defined as
\begin{equation}
  A(v)=-\Tr\left(\frac{\partial \E(v^m)}{\partial
y}\right),
\label{def-opearor}\end{equation}
with domain
\begin{equation}
  D(A)=\{v\in L^1(\Gamma)\cap L^\infty(\Gamma)\,:\,A(v)\in L^1(\Gamma),\,\|v\|_{L^\infty(\Gamma)}\le
  \|f\|_{L^\infty(\Gamma)}\}
\end{equation}
This operator is nothing but the half-laplacian of the power $m$,
$A(v)=(-\Delta)^{1/2}v^m$.

\subsection{The elliptic problem}\label{sect-elliptic}

Therefore, in order to perform the plan we need to establish the solvability and properties
of the elliptic problem
\begin{equation}
\left\{
\begin{array}{ll}
\Delta w=0\qquad &\mbox{ in } \Omega,\\
-\dfrac{\partial w}{\partial y}+w^{1/m}=g\qquad&\mbox{ on } \Gamma.
\end{array}
\right. \label{pe:local}
\end{equation}
for all $g\in L^1(\Gamma)\cap L^\infty(\Gamma)$. As we have said before, the
power $w^{1/m}$ in the boundary condition means $|w|^{1/m-1}w$ if $w$ takes on
some negative values. We will also prove that if $g\ge0$ then $w\ge0$. A weak
solution to this problem is a function $w\in W_{loc}^{1,1}(\Omega)$, such that
$\Tr(w^{1/m})\in L^1(\Gamma)$, verifying
\begin{equation}
\int_{\Omega}\langle\nabla
w,\nabla\varphi\rangle+\int_{\Gamma}w^{1/m}\varphi-\int_{\Gamma}g\varphi=0
\label{eq:weak-elliptic}
\end{equation}
for any $\varphi\in C_0^1(\overline\Omega)$. We have to prove existence of the
solution $w$ and contractivity of the map $g\mapsto \Tr(w^{1/m})$ in the norm of
$L^1(\Gamma)$, which plays the role of $\cal X$ in the definition of
accretivity.

To prove this we  perform an approximation substituting the unbounded domain
$\Omega$ by an increasing sequence of bounded domains $\Omega_R$ (half balls),
imposing zero Dirichlet condition on the part of the boundary of the domain
which does not lie on the hyperplane $y=0$. The approximate problems are
\begin{equation}
\left\{
\begin{array}{ll}
\Delta w=0\qquad &\mbox{ in } \Omega_R=\Omega\cap B_R,\\
\dfrac{\partial w}{\partial
y}=w^{1/m}-g\qquad&\mbox{ on } \Gamma_R=\partial\Omega_R\cap\{y=0\},\\
w=0\qquad&\mbox{ on } \Sigma_R=\partial\Omega_R\cap\{y>0\}\,,
\end{array}
\right. \label{pe:local.bola}
\end{equation}
where $B_R=B_R(0)$. The concept of weak solution for a given datum $g\in L^1(\Gamma_R)\cap
L^\infty(\Gamma_R)$ is analogous to the one given above
\eqref{eq:weak-elliptic}, after changing the domains of the
integrals into the corresponding bounded domains.

\begin{Theorem}\label{th:pe}
For every $g\in L^1(\Gamma)\cap L^\infty(\Gamma)$ there exists a unique weak
solution $w\in H^1(\Omega)$ to Problem \eqref{pe:local} such that
$\Tr(w^{1/m})\in L^1(\Gamma)\cap L^\infty(\Gamma)$. Moreover, if $w$ and
$\widetilde w$ are the solutions corresponding to data $g$ and $\widetilde g$,
then
\begin{equation}
\int_{\Gamma}\left[w^{1/m}(x,0)-\widetilde
w^{1/m}(x,0)\right]_+\,dx\le
\int_{\Gamma}\left[g(x)-\widetilde
g(x)\right]_+\,dx.
\label{eq:contract-e}\end{equation}
This in turn implies that if $g\ge0$ in $\Gamma$ then $w\ge0$ in
$\overline\Omega$. Moreover, $\|u\|_{L^\infty(\Gamma)}\le
  \|g\|_{L^\infty(\Gamma)}$.
\end{Theorem}

\noindent{\it Proof. }
\noindent {\sc Step 1}. We first prove that there
exists a weak solution $w\in H^1(\Omega_R)$ to Problem
\eqref{pe:local.bola}. This is done by solving the following minimization
problem:

\begin{center}\begin{minipage}{14cm}To find a function $w\in H^1(\Omega_R)$ minimizing
$$
J(w)=\frac12\int_{\Omega_R}|\nabla
w|^2+\frac{m}{m+1}\int_{\Gamma_R}|w|^{\frac{m+1}m}-
\int_{\Gamma_R}wg.
$$
\end{minipage}\end{center}

This functional is coercive, since
$$
J(w)\ge C_1 \| w\|^2_{H^1(\Omega_R)}-C_2 \| w\|_{H^1(\Omega_R)},
$$
which follows by using the Poincar{\'e} inequality, Cauchy-Schwartz and the trace
embedding. Moreover, coercivity then provides a  bound for $\|
w\|_{H^1(\Omega_R)}$, though it depends on $R$.

{\sc Step 2.} We now establish contractivity of solutions to Problem
\eqref{pe:local.bola} in $L^1(\Gamma_R)$.  Let
$w$ and $\widetilde w$ be the solutions corresponding to data $g$ and
$\widetilde g$. We claim that
\begin{equation}
\int_{\Gamma_R}\left[w^{1/m}(x,0)-\widetilde
w^{1/m}(x,0)\right]_+\,dx\le
\int_{\Gamma_R}\left[g(x)-\widetilde
g(x)\right]_+\,dx. \label{eq:contract-r}\end{equation} This inequality follows
easily if we consider in the weak formulation the test function
$\varphi=p(w-\widetilde w)$, where $p$ is any smooth monotone approximation of
the sign function, $0\le p(s)\le 1$, $p'(s)\ge 0$. We get
$$
\int_{\Omega_R}p'(w-\widetilde
w)|\nabla(w-\widetilde w)|^2+\int_{\Gamma_R}(w^{1/m}-\widetilde
w^{1/m})\,p(w-\widetilde w)-\int_{\Gamma_R}(g-\widetilde g)\,p(w-\widetilde
w)=0.
$$
Passing to the limit, we obtain
$$
\int_{\Gamma_R}(w^{1/m}-\widetilde
w^{1/m})_+\,dx\le\int_{\Gamma_R}(g-\widetilde g)
\,\mathop{\rm sign}(w-\widetilde
w)\,dx\le\int_{\Gamma_R}(g-\widetilde g)_+\,dx\,.
$$
In particular, under the assumption $g\ge0$ we have $w(\cdot,0)\ge0$. Moreover,
$w(\cdot,0)\in L^1(\Gamma_R)\cap L^\infty(\Gamma_R)$. Finally, since the
Poisson kernel in the half-ball is nonnegative, we also conclude that $w\ge0$
in $\Omega_R$.

{\sc Step 3.} In order to pass to the limit $R\to\infty$  in the
case of nonnegative data, we use a monotonicity property of the
family of approximate solutions, denoted here by~$w_R$. Namely,
$R<R'$ implies $w_R\le w_{R'}$ in $\Omega_{R'}$. This follows from
the ordering of the restrictions,  using again that the Poisson
kernel is nonnegative. The ordering of the restrictions results from
comparison in $\Gamma_{R}$, since $w_{R'}\ge0$ in $\Sigma_R$ (there
is a contraction property analogous to \eqref{eq:contract-r} for
problems with non-homogeneous boundary data).

Monotonicity implies that there exists the pointwise (and also in
the sense of distributions) limit $w=\lim\limits_{R\to\infty}w_R$.
This limit satisfies $w\ge0$ in $\Omega$, $w^{1/m}(\cdot,0)\in
L^1(\Gamma)\cap L^\infty(\Gamma)$. Since $|\nabla w_R|$ is uniformly
bounded in $L^2(\Omega_R)$,
$$
\int_{\Omega_R}|\nabla w_R|^2\le \int_{\Gamma_R} gw_R\le
\|g\|_{L^1(\mathbb{R}^N)}\|g\|_{L^\infty(\mathbb{R}^N)}^m\,,
$$
we conclude that $\nabla w_R\rightharpoonup \nabla w$ in $L^2(\Omega)$. This is
enough to pass to the limit in the identity
$$\int_{\Omega_R}\langle\nabla
w_R,\nabla\varphi\rangle+\int_{\Gamma_R}w_R^{1/m}\varphi-\int_{\Gamma_R}g\varphi=0,
$$
to show that $w$ satisfies \eqref{eq:weak-elliptic}. Also the estimate of the
$L^2$ norm of the gradients passes to the limit, and leads to
\begin{equation}
\int_{\Omega}|\nabla
w|^2\le\|g\|_{L^1(\mathbb{R}^N)}\|g\|_{L^\infty(\mathbb{R}^N)}^m\,.
\label{rabia}\end{equation}
{\sc Step 4.} The pass to the limit in the case of non-positive data uses a
similar argument. Finally, in the case of data $g$ of both signs, we use
comparison with the solutions with data $g_1=g^+\ge 0$ and $g_2=-g^-\le 0$ and
compactness to pass to the limit.

{\sc Step 5.} Contractivity for the limit problem is proved exactly in the same
way as for the approximate problems. This gives uniqueness. We also have that
the $L^1$ norm and the $L^\infty$ norm of the function $w^{1/m}(\cdot,0)$ are
bounded respectively by the $L^1$ norm and the $L^\infty$ norm of the datum.
\qed


\subsection{Existence of solution for the evolution problem}\label{sect-exist}

We now use the previously mentioned procedure to construct the solution to the
evolution problem \eqref{pp:local}. We recall that the Crandall-Liggett result
only provides us in principle with an abstract type of solution called {\sl mild solution.}

\begin{Theorem}\label{th:plocal}
For every $f\in L^1(\Gamma)\cap L^\infty(\Gamma)$ there exists a weak solution
$(u,w)$ to Problem \eqref{pp:local} with $u(\cdot,t)\in L^1(\Gamma)\cap
L^\infty(\Gamma)$ for every $t>0$ and $w\in L^2([0,T];H^1(\Omega))$. Moreover,
the following contractivity property holds: if $(u,w),\,(\widetilde
u,\widetilde w)$ are the constructed weak solutions corresponding to initial
data $f,\,\widetilde f$, then
\begin{equation}
\int_{\Gamma}[u(x,t)-\widetilde u(x,t)]_+\,dx\le
\int_{\Gamma}[f(x)-\widetilde f(x)]_+\,dx.
\label{eq:contract-p}\end{equation}
In particular a comparison principle for constructed solutions is obtained.
\end{Theorem}

\noindent{\it Proof. }
For each time $T>0$ we divide the time interval $[0,T]$ in $n$ subintervals.
Letting $\varepsilon=T/n$, we construct the function $w_\varepsilon$ piecewise
constant in each interval $(t_{k-1},t_k]$, where $t_k=k\varepsilon$,
$k=1,\cdots,n$, as the solutions to the discretized Problems
\eqref{pe:local-iterate}. For convenience we write here again the problems:
$w_{\varepsilon,k}$ solves
$$
\left\{
\begin{array}{ll}
\Delta w_{\varepsilon,k}=0\qquad &\mbox{ in } \Omega,\\
\varepsilon \dfrac{\partial w_{\varepsilon,k}}{\partial
y}=u_{\varepsilon,k}- u_{\varepsilon,k-1}\qquad&\mbox{ on }
\Gamma,
\end{array}
\right.
$$
with $u_{\varepsilon,0}=f$. Our solution is the (uniform in $[0,T]$) limit
$$
(u,w)=\lim_{\varepsilon\to0}(u_\varepsilon,w_\varepsilon)
$$
in $L^1_{loc}(\Omega)$. It is a {\sl mild solution}, whose existence is
guaranteed by the classical semigroup approach, $u\in C([0,\infty):
L^1(\Gamma))$. In fact we obtain first the function $u$ by  Crandall-Ligget's
Theorem, and the harmonic extension $w$ of $u^m$ coincides with
$\lim\limits_{\varepsilon\to0}w_\varepsilon$. By construction we have $w\in
L^\infty(\Omega\times[0,T])$. We must now show that we have obtained in fact a
weak solution.

Multiplying the equation by $w_{\varepsilon,k}$, integrating by
parts, and applying Young's inequality, we obtain
\begin{equation}
\label{eq:energy.epsilon}
\varepsilon\int_\Omega |\nabla
w_{\varepsilon,k}|^2\,d\overline x\le \dfrac1{(m+1)}\Big(\int_\Gamma
|u_{\varepsilon,k-1}|^{m+1}\,dx-\int_\Gamma |u_{\varepsilon,k}|^{m+1}\,dx\Big).
\end{equation}
Adding from $k=1$ to $k=n$ we get
$$
\int_0^T\int_\Omega |\nabla w_\varepsilon(\overline x,t)|^2\,d\overline xdt\le
\dfrac1{(m+1)}\int_\Gamma |f(x)|^{m+1}\,dx.
$$
Passing to the limit, the same estimate is obtained for $|\nabla w|$, and
therefore $w\in L^2([0,T];H^1(\Omega))$. On the other hand,
\eqref{eq:energy.epsilon} yields
$$
\displaystyle \int_\Gamma |u_{\varepsilon,k}|^{m+1}\,dx\le
\int_\Gamma |u_{\varepsilon,k-1}|^{m+1}\,dx \le \int_\Gamma
|f|^{m+1}\,dx,
$$
which, after passing to the limit, gives
$$
\int_\Gamma |u(x,t)|^{m+1}\,dx\le  \int_\Gamma |f(x)|^{m+1}\,dx,
$$
for every $t\in[0,T]$. Now, choosing appropriate test functions, as in
\cite{RV}, it follows that we can pass to the limit in the elliptic weak
formulation to get the identity of the parabolic weak formulation. We first
have
$$
\int_\Omega\langle\nabla
w_{\varepsilon,k},\nabla\varphi\rangle=
\int_\Gamma\frac1\varepsilon(u_{\varepsilon,k-1}-u_{\varepsilon,k})\varphi.
$$
Integrating in $(t_{k-1},t_k)$ and adding for $k=1,\cdots,n$, we get that the
right-hand side becomes
$$
\begin{array}{l}
\displaystyle\sum_{k=1}^n\int_{t_{k-1}}^{t_k}\int_\Gamma\frac1\varepsilon(u_{\varepsilon,k-1}(x,t)-
u_{\varepsilon,k}(x,t))\varphi(x,0,t)\,dxdt= \\ [4mm] \qquad\displaystyle=
\int_0^T\int_\Gamma
u_\varepsilon(x,t)\,\frac1\varepsilon(\varphi(x,0,t+\varepsilon)-\varphi(x,0,t))\,dxdt\\
[4mm] \qquad\displaystyle+
\frac1\varepsilon\int_0^\varepsilon\int_\Gamma f(x)\varphi(x,0,t)\,dxdt-
\frac1\varepsilon\int_{T-\varepsilon}^T\int_\Gamma u_\varepsilon(x,T)\varphi(x,0,t)\,dxdt.
\end{array}
$$
Passing to the limit $\varepsilon\to0$ we get \eqref{weak-local}.

\

The contractivity
\eqref{eq:contract-p} obtained in Theorem~\ref{th:pe} in each step is inherited in the
limit. In fact, if $w_{\varepsilon,k}$ and $\widetilde w_{\varepsilon,k}$ are
the discretized approximations of $w$ and $\widetilde  w$, then we have
$$
\int_{\Gamma}\left[u_{\varepsilon,k+1}(x)-\widetilde u_{\varepsilon,k+1}(x)\right]_+\,dx\le
\int_{\Gamma}\left[u_{\varepsilon,k}(x)-\widetilde
u_{\varepsilon,k}(x)\right]_+\,dx,
$$
which easily implies \eqref{eq:contract-p}. Comparison is a trivial consequence
of contractivity. \qed

\noindent{\it Remark. }
This contractivity property also implies the following estimates for the weak
solution to Problem \eqref{eq:main} just constructed
\begin{equation}
\|u(\cdot,t)\|_{L^1(\mathbb{R}^N)}\le\|f\|_{L^1(\mathbb{R}^N)}\,,\qquad
\|u(\cdot,t)\|_{L^\infty(\mathbb{R}^N)}\le\|f\|_{L^\infty(\mathbb{R}^N)}\,.
\label{norms}\end{equation}
Using now the Poisson kernel of the half-space and Young's inequality, we have
that for every $y>0$ and every $1/m\le p\le\infty$, it holds
\begin{equation}
  \|\E(u^{m})(\cdot,y,t)\|_{L^p(\mathbb{R}^N)}\le\|u^{m}(\cdot,0,t))\|_{L^p(\mathbb{R}^N)}\le
  \|f\|_{L^\infty(\mathbb{R}^N)}^{m-1/p}\|f\|_{L^1(\mathbb{R}^N)}^{1/p}\,,
\label{eq:norm-ext}\end{equation}
with the first inequality replaced by equality if $p=1$, $m>1$.

We end this section with a property satisfied by the weak solutions just
constructed which is very useful in the proofs to come, with a number of other
applications.

\begin{Proposition}\label{pro:ut+u} Assume $f\ge0$. If $w$ is the nonnegative weak solution to Problem
\eqref{pp:local} constructed in Theorem~{\rm \ref{th:plocal}}, then the inequality
\begin{equation}
(m-1)t\frac{\partial w}{\partial t}+ m\,w \ge 0
\label{wt+w>}\end{equation} holds in the sense
of distributions.
\end{Proposition}
\noindent{\it Proof.} We use a simple scaling argument based on the
homogeneity of the problem, (see \cite{BenilanCrandall:1}). Assume
first $m>1$. We have that for every $\lambda>1$ the function
$w_\lambda(x,t)=\lambda w(x,\lambda^{\frac{m-1}m} t)$ is also a
solution with initial value bigger than $w(\cdot,0)$. Then by the
comparison principle we get,
$$
\lim\limits_{\lambda\to1^+}\dfrac{w_\lambda(x,t)-w(x,t)}{\lambda-1}\ge0,
$$
which gives
\eqref{wt+w>}. For $m<1$ the sign is  reversed. Comparison can be easily justified in the discretized approximations.
Since the function $w_\lambda$ is the limit of the rescaled approximations of
the function $w$, we can compare $w$ and $w_\lambda$.
\qed

We observe for future reference that at the nonlocal level of function $u$ we
have the ``monotonicity formulae''
\begin{equation}
\begin{array}{ll}
\dfrac{\partial u}{\partial t}\ge -\dfrac{u}{(m-1)t}&\quad\mbox{if } m>1, \\
[3mm]
\dfrac{\partial u}{\partial t}\le \dfrac{u}{(1-m)t}&\quad\mbox{if } m<1.
\end{array}
\label{ut+u>}\end{equation}
Formula \eqref{wt+w>} is empty for $m=1$, but in this case it is easy to derive
from the explicit representation of the solution that $t\partial_tu+Nu\ge0$.

\noindent{\it Remark. } In the PME model (the local analogue), a similar lower estimate of
$\partial_tu$ is also available in the case $m<1$. The proof uses in an
essential way a second variable called the pressure, which is a potential for
the velocity. It is not clear which could be the corresponding pressure for the
nonlocal problem.

\section{Uniqueness of weak energy solutions}\label{sect-uniq}
\setcounter{equation}{0}

In the previous section we have constructed a weak solution to the
local Problem~\eqref{pp:local}. As we have said, the construction
itself shows that this weak solution is in fact a weak energy
solution. We prove uniqueness using an argument taken from  Oleinik
et al.~\cite{OKC}.

\begin{Lemma}
 Assume $f\in L^1(\mathbb{R}^N)\cap L^\infty(\mathbb{R}^N)$.  There is at most one weak energy solution to Problem
  \eqref{eq:main}.
\label{oleinik}\end{Lemma}

\noindent{\it Proof. } Let $(u,w)$ and $(\widetilde u,\widetilde w)$
be two weak solutions to Problem \eqref{pp:local}. We take as test,
in the weak formulation, the following function
$$
\varphi(\overline x,t) =\int_t^T (w-\widetilde w)(\overline x,s)\,ds,\qquad 0\le t\le
T,
$$
with $\varphi\equiv0$ for $t\ge T$. Observe that this is a good test function
when $w$ and $\widetilde w$ are weak energy solutions. We have
$$
\begin{array}{l}
\displaystyle\int_0^T\int_{\Omega}\langle\nabla(w-\widetilde w)(\overline x,t),\int_t^T\nabla(w-\widetilde
w)(\overline x,s)\,ds\rangle \,d\overline xdt \\ [4mm]
+\displaystyle\int_0^T\int_{\Gamma}(u-\widetilde  u)(x,t)(u^m-\widetilde  u^m)(x,t)\,dxdt = 0.
\end{array}
$$
Integration of the first term gives
$$
\begin{array}{l}
\displaystyle\frac 12\int_{\Omega}\Big|\int_0^T \nabla(w-\widetilde  w)(\overline x,s)\, ds\Big|^2d\overline x \\
[4mm] +
\displaystyle\int_0^T\int_{\Gamma}(u-\widetilde  u)(x,t)(u^m-\widetilde  u^m)(x,t)\,dxdt = 0.
\end{array}$$
Since both terms are nonnegative, they must be zero. Therefore, $u=\widetilde u$
in $\Gamma$. Obviously this also gives, as a byproduct, $w=\widetilde  w$ in
$\Omega$. \qed

\noindent {\it Remark. }   Observe that this proof only requires $u(\cdot,t),\,\widetilde u(\cdot,t)\in
  L^1(\mathbb{R}^N)\cap L^{m+1}(\mathbb{R}^N)$.

\section{Properties of weak energy solutions}\label{sect-props}

In the next sections we restrict ourselves to the range  $m>m_*$.

\subsection{Conservation of mass}

We establish next a property that is typical of diffusive processes.

\begin{Theorem}\label{th:mass}
If $(u,w)$ is a weak energy solution to Problem \eqref{pp:local}
with initial datum $f$, then for every $t>0$ we have
\begin{equation}
\int_{\mathbb{R}^N}u(x,t)\,dx=
\int_{\mathbb{R}^N}f(x)\,dx.
\label{eq:cons-mass}\end{equation}
\end{Theorem}

\noindent{\it Proof.} The case $m=1$ follows from the explicit representation \eqref{linear-rep}.
For general $m$ we use the integral identity
\eqref{weak-local} with a particular test function.
Consider a nonnegative nonincreasing cut-off function $\psi(s)$ such that
$\psi(s)=1$ for $0\le s\le1$, $\psi(s)=0$ for $s\ge2$, and define
$\varphi_R(\overline x)=\psi(|\overline x|/R)$. We obtain, for every
$t_2>t_1\ge0$,
$$
\int_\Gamma \Big(u(x,t_2)-u(x,t_1)\Big)\varphi_R(x,0)\,dx=
-\int_{t_1}^{t_2}\int_\Omega\langle\nabla w(\overline
x,t),\nabla\varphi_R(\overline x)\rangle\,d\overline xdt.
$$
Integrating by parts, noting that $\dfrac{\partial\varphi_R}{\partial
y}(x,0)=0$, we have that the space integral inside the right-hand side is
$$
I=\int_\Omega w(\overline x,t)\Delta\varphi_R(\overline x)\,d\overline x.
$$
With the test function chosen we have
$$
|I|\le cR^{-2}\int_{R<|\overline x|<2R}w(\overline x,t)\,d\overline x.
$$
We now estimate this integral for every $t_1<t<t_2$. If $m>1$ we just observe
that
$$
\begin{array}{rl}
|I|&\le\displaystyle cR^{-2}\int_0^{2R}\int_{\mathbb{R}^N}|w(x,y,t)|\,dxdy\le
cR^{-1}\|u(\cdot,t)\|^{m-1}_{L^\infty(\mathbb{R}^N)}\|u(\cdot,t)\|_{L^1(\mathbb{R}^N)} \\
[3mm] &\le cR^{-1}\to0\,,
\end{array}
$$
where we have used \eqref{eq:norm-ext}. In the case  $m<1$, applying H\"older's inequality with some exponent $p>1/m$ we have
$$
\begin{array}{rl}
|I|&\displaystyle\le
cR^{-2}\|u(\cdot,t)\|_{L^\infty(\mathbb{R}^N)}^{m-1/p}\Big(\int_{R<|\overline
x|<2R}|w(\overline x,t)|^{1/m}\,d\overline x\Big)^{1/p}
\,|\{R<|\overline x|<2R\}|^{(p-1)/p} \\ [3mm] &\displaystyle\le
cR^{-2+(N+1)((p-1)/p)}\|u(\cdot,t)\|_{L^\infty(\mathbb{R}^N)}^{m-1/p}\Big(\int_0^{2R}\int_{\mathbb{R}^N}|w(x,y,t)|^{1/m}\,dxdy\Big)^{1/p}
\\ [3mm] &\displaystyle\le
cR^{-2+(N+1)((p-1)/p)+1/p}\|u(\cdot,t)\|_{L^\infty(\mathbb{R}^N)}^{m-1/p}\|u(\cdot,t)\|_{L^1(\mathbb{R}^N)}^{1/p}\le
cR^{N(p-1)/p-1}.
\end{array}
$$
Finally observe that if $m>m_*$ we can choose $1/m<p<N/(N-1)$ to force the last
term to go to zero as $R\to\infty$.
\qed

\noindent {\sc Solutions that lose mass.} We now present a result that shows the
necessity of the condition $m>m_*$. In fact if $0<m<m_*$ there is a phenomenon
of extinction in finite time, which makes impossible to have conservation of
mass. The proof is almost exactly the same as the one in
\cite{Benilan-Crandall} for the  PME model, where the
corresponding condition on $m$ is $0<m<(N-2)/N$ instead of $0<m<m_*=(N-1)/N$.

\begin{Proposition}
  Let $N>1$ and $0<m<(N-1)/N$, and let $f\in L^{1}(\mathbb{R}^N)\cap L^{(1-m)N}(\mathbb{R}^N)$. Then
  there is a finite time $T>0$ such that the solution $u$ to Problem
  \eqref{eq:main} satisfies $u(x,T)\equiv0$ in $\mathbb{R}^N$.
\label{pro:extinction}\end{Proposition}

\noindent{\it Proof.}
As we have said, the proof follows the argument in \cite{Benilan-Crandall}.
Therefore we leave the details to be consulted there. Assume also for
simplicity $u\ge0$.

If we consider $\varphi=w^p$ in the equality \eqref{strong-local}, we get
$$
\frac{4p}{(p+1)^2}\int_\Omega |\nabla w^{\frac{p+1}2}|^2\,d\overline x+\frac1{1+pm}\int_\Gamma
\dfrac{\partial(u^{pm+1})}{\partial t}\,dx=0.
$$
The use of this test function is justified in~\cite{Benilan-Crandall}. By the
trace inequality we obtain that there exists a positive constant $C=C(p,N)$
such that
$$
C\Big(\int_\Gamma u^{\frac{(p+1)mN}{N-1}}\,dx\Big)^{\frac{N-1}N}+\int_\Gamma
\dfrac{\partial(u^{pm+1})}{\partial t}\,dx\le 0.
$$
If we now choose $p=(N(1-m)-1)/m>0$, we get $((p+1)mN)/(N-1)=pm+1$.
Therefore, the function $J(t)=\int_\Gamma u^{pm+1}\,dx$ satisfies the
inequality
$$
J'(t)+CJ^{\frac{N-1}N}(t)\le0.
$$
This implies extinction in finite time for $J(t)$ and thus for $u(x,t)$,
provided $J(0)$ is finite.
\qed

\noindent  {\sc Example.} There exists an explicit example of the above extinction
property for a particular $m<m_*$. It has the form (separated variables),
$$
u(x,t)=G(x)H(t).
$$
Substituting this expression in \eqref{eq:main}, we have
$H(t)=c(T-t)^{1/(1-m)}$, and $G$ solves the nonlocal equation
$$
(-\Delta)^{1/2} G^m=G.
$$
In the special case $m=(N-1)/(N+1)$, there exists an explicit family of
solutions
$$
G_{\delta,\tau}=A(\tau)[\tau^2+|x-\delta|^2]^{-(N+1)/2},
$$
with $\tau>0$, $\delta\in\mathbb{R}^N$ and where $A(\tau)>0$ has an explicit
expression, see
\cite{BLW} and also
\cite{CLO}. Observe that $m<m_*$ and $G\in L^1(\mathbb{R}^N)$, and
thus $u(\cdot, t)\in L^1(\mathbb{R}^N)$ for any $0\le t<T$, while
$u(x,T)\equiv0$.

\subsection{Positivity and regularity}

We show in this subsection that nonnegative bounded weak energy solutions are
in fact positive everywhere in $\overline{\Omega}$. This is true for every
$m>m_*$. The result is in sharp contrast with what happens for the local analog, the
 PME $\partial_t u=\Delta u^m$ in the case
$m>1$, for which initial values with compact support produce solutions that
develop a free boundary. Free boundaries are a main feature of the standard
 PME theory, but they are not available here.

\medskip

The idea behind our positivity result is as follows: if $u\ge0$ is a
classical solution and $u(x_0,t)=0$ for some $x_0$ and $t$, then
formula \eqref{def-riesz} gives $(-\Delta)^{1/2}u^m(x_0,t)<0$, and
hence $\partial_t u(x_0,t)>0$. If $m>1$ we only have that the
bounded solution $u$ is H\"{o}lder continuous by \cite{AC}. In the case $m<1$ we do not even have that regularity.
Now we perform rigourously the proof of positivity for weak
solutions using the extension Problem~\eqref{pp:local}.

\begin{Theorem}\label{th:positivity}
Assume $m>m_*$. The weak energy solution $u$ to Problem
\eqref{eq:main} with a bounded nonnegative initial datum is positive
for positive times. Even more, the corresponding pair $(u,w)$ to
Problem \eqref{pp:local} satisfies $w(x,y,t)>0$ for every
$x\in\mathbb{R}^N$, $y\ge0$ and every $t>0$.
\end{Theorem}

\noindent{\it Proof. } The case $m=1$ follows from formula \eqref{linear-rep}.

{\sc Case $m>1$}: We already know that for every $t>0$, function $w$ is
H\"{o}lder continuous, $w\ge 0$ in $\overline\Omega$, and $w>0$ in $\Omega$. We
have to prove that $w>0$ also on the boundary, $\Gamma$. By the  comparison result,  we only
need to consider compactly supported initial data.

In a first step we show that if $w$ is not strictly positive on $\Gamma$ at
some time $T$, then the supports of $w(\cdot,0,t)$ form an expanding (in time)
family of compact sets for $0\le t\le T$. This follows from estimate
\eqref{ut+u>} and Alexandrov's reflection principle. In fact \eqref{ut+u>}
implies that
$$
w(\overline x,t_2)\ge w(\overline x,t_1)\Big(\frac{t_1}{t_2}\Big)^{m/(m-1)},
$$
for every $\overline x\in\overline\Omega$, $t_2\ge t_1>0$, which is called
retention property. Next we claim that if the support of the initial value $f$
is contained in the ball $B_R$, and $w(x_0,0,t_0)=0$ for some
$x_0\in\mathbb{R}^N$, $t_0>0$, then the support of $w(\cdot,0,t)$ is also
compact for any $0<t\le t_0$, and contained in a ball of radius depending on
$|x_0|$ and $R$.

To prove the claim we reflect, for any given point $(x_1,y)\in\overline\Omega$,
around the hyperplane $\pi\equiv(x_1-x_0)\cdot(x-(x_0+x_1)/2)=0$. It is clear
that if $|x_0-(x_0+x_1)/2|>|x_0|+R$, then the hyperplane $\pi$ divides the
half-space $\overline\Omega$ in two parts, $\overline\Omega=M_1\cup M_2$ with
$B_R\times[0,\infty)\subset M_1$, $(x_0,0)\in M_2$. In this way, by the
comparison principle, we obtain that the function
$z(x,y,t)=w(x,y,t)-w(x_0+x_1-x,y,t)$ satisfies $z(x,y,t)\ge0$ for every
$(x,y)\in M_1$, $t>0$. The comparison principle holds on the half-space $M_1$.
Thus $w(x_1,0,t)\le w(x_0,0,t)=0$ for every $0< t\le t_0$. A sufficient
condition for $x_1$ to get this argument work is $|x_1|\ge3|x_0|+2R$.

In a second step we assume (thanks to the previous argument) that
$w(x,0,t)\equiv0$ in some ball $B\subset\Gamma$ for $0\le t\le t_1$. Then we
have, for every test function $\varphi$ that vanishes on $\partial
B\times(0,\infty)$, that
$$
0=\int_0^{t_1}\int_B u\frac{\partial\varphi}{\partial t}\,dxdt=
\int_0^{t_1}\int_0^\infty\int_B\langle\nabla w,\nabla\varphi\rangle\,dxdydt=
-\int_0^{t_1}\int_B\frac{\partial w}{\partial y}\,\varphi\,dxdt.
$$
This gives $\partial_y w(x,0,t)\equiv0$ for $x\in B$, $0\le t\le t_1$.  But $w$
is a continuous nonnegative harmonic function in the half cylinder which
vanishes on the part of the boundary $y=0$. Hopf's Lemma implies $\partial_y
w(x,0,t)<0$ for $x\in B$. This is a contradiction. Therefore, $w(\cdot,0,t)$ is
positive everywhere.

{\sc Case $m<1$}: In this case the proof is different, based on a weak Harnack
inequality. Using estimate \eqref{ut+u>} and the fact that the solution is
bounded, we know that for every $t>0$ there exists a constant $A>0$ such that
$$
\int_\Omega \langle\nabla w,\nabla\varphi\rangle+A\int_\Gamma w\varphi\ge0,
$$
for every nonnegative test function $\varphi$. Once we have this, we can use part of the
proof of Lemma~2.4 in \cite{Cabre-Sola} to get a weak Harnack inequality in
each large ball $B_R=\{|x-x_0|^2+y^2<R^2\}\subset\mathbb{R}^{N+1}$ with center
on $\Gamma$. First of all, conservation of mass plus
\eqref{eq:norm-ext}, together with  the fact that $w\ge0$, $w\not\equiv0$, imply
that there exists some $R>0$ large such that $\int_{B_{R/2}}w\,d\overline x>0$.

If we consider the function
$$
z(x,y)=e^{-A|y|}w(x,|y|),
$$
then we have that $z$ satisfies
$$
\int_{B_R}\langle\nabla z,\nabla\varphi\rangle-2A\int_{B_R}\mbox{sign}(y)\frac{\partial
z}{\partial y}\varphi-A^2\int_{B_R}z\varphi\ge0,
$$
i.e., it is a weak supersolution to an equation for which we can apply
Theorem~8.18 in
\cite{GT}, to get
$$
\inf_{B_{R/2}}z\ge cR^{-N-1}\|z\|_{L^1(B_{R/2})}.
$$
This means $z>0$, (and thus $w>0$) in $\overline\Omega$.
\qed

As a corollary of this result, we can establish also regularity for the case
$m<1$, provided $m>m_*$.

\begin{Theorem}\label{th:regularity2}
Let $m_*<m<1$. Then any bounded weak energy nonnegative solution $u$ to Problem
\eqref{eq:main} satisfies $u\in C^\alpha(\mathbb{R}^N\times(0,\infty))$ for some $0<\alpha<1$.
\end{Theorem}

\noindent{\it Proof. }
The above-mentioned regularity result of \cite{AC}  applies to the equation
$$
\dfrac{\partial \beta(v)}{\partial t} + (-\Delta)^{1/2} (v)=0
$$
in some ball $x\in B\subset\mathbb{R}^N$, $t>0$, with some nondegeneracy
condition on the constitutive monotone function $\beta$. Once we know that in
such a ball the solution is bounded below away from zero, the requirements on
the function $\beta$ in \cite{AC} are fulfilled.
\qed

\subsection{Strong solutions}\label{sect-smooth}

We prove here that every nonnegative bounded weak energy solution is in fact a
strong solution. We need to show that the time partial derivative of $u$ is an
$L^1$ function and that the second equality in \eqref{pp:local} holds almost
everywhere.

As a first step we show that the time-increment quotients are bounded in
$L^1(\Gamma)$, and thus the limit $\partial_t u$ must be a Radon measure. Our
purpose is to prove that the limit is still in $L^1$, and this is proved later.
Observe that the result is clear in the case $m=1$, from \eqref{linear-rep}. We
therefore assume $m\neq1$.

\begin{Proposition}
  If $u$ is the weak solution to Problem \eqref{eq:main} constructed in Theorem~{\rm \ref{th:main}}, then
  $h^{-1}(u(\cdot,t+h)-u(\cdot,t))\in L^1(\mathbb{R}^N)$ for every $t,\,h>0$.
\end{Proposition}

The proof is exactly the same as in the PME case, see
\cite{BenilanCrandall:1}. In fact, the following estimate holds
\begin{equation}
\int_{\mathbb{R}^N}\dfrac1h\Big|u(x,t+h)-u(x,t)\Big|\,dx\le\frac{2}{|m-1|t}\,\|f\|_{L^1(\mathbb{R}^N)}.
\label{utmeasure}
\end{equation}
An analogous result for more general nonlinearities is given in
\cite{CP}.

Observe that by the Mean Value Theorem, the same type of estimate can be
obtained for the time-increment quotients of $w=u^m$ if $m>1$ since $u\in
L^\infty$. If $m<1$ we obtain a bound locally for nonnegative solutions, since
in each ball of $\mathbb{R}^N$ we have $u\ge c>0$.

We now prove a result which turns out to be fundamental in the proof
of Theorem~\ref{th:main2}.
\begin{Proposition}
Let $f\ge0$ and let $u$ be the weak solution to Problem \eqref{eq:main}
constructed in Theorem~\ref{th:main}. Then $\partial_t(u^{(m+1)/2})\in
L^2_{loc}(\mathbb{R}^N\times[0,T])$.
\label{pro:estimates}\end{Proposition}

\noindent{\it Proof. } The formal proof is simple, using the function $\partial_t w$
as test function in the weak formulation, much as in the PME case,
cf.~\cite[Section 5.5]{Vazquez:Book}, of course after some obvious changes. The
problem is the justification of the calculations, since we have not established
the existence of any kind of differentiability for the solutions or suitable
approximations.

Here we do by brute force as follows:

$\bullet$ We use the weak formulation with the test function $\varphi=\delta^h(w*\rho)$,
where $\rho=\eta*\eta$ and $\eta$ is a convolution kernel, acting only on the time variable.
We fix the notation $\hat f=f*\eta$, $\widetilde f=f*\rho$, applied to functions of $t$. We also make use of the
following calculus identity
\begin{equation}
\int\widetilde f(t)g(t)\,dt=\int\hat f(t)\hat g(t)\,dt.
\label{convol}\end{equation} We prove this identity in
Lemma~\ref{lem-a1} of the  Appendix for the reader's convenience. In
addition, we take $\eta(t)=\eta_h(t)=(1/h)\eta_1(t/h)$, where
$\eta_1$ is a smooth, symmetric, compactly supported nonnegative
function, with support $[-1,1]$. Thus, $\rho=\rho_h$ inherits the
same properties, with supp$(\rho_1)=[-2,2]$. We also use the
notation
$$
\delta^h w(t)=\dfrac1{2h}\Big(w(t+h)-w(t-h)\Big)
$$
for a discrete time derivative, omitting the spatial variable. Observe that
since $|\nabla w|$ is in $L^2$, this is a good test function.

$\bullet$  Going  to the weak formulation
\eqref{weak-local2} with the  test function $\delta^h\widetilde
w,$ we have
\begin{equation}
-\int_{t_1}^{t_2}\int_\Gamma\delta^h(w*\rho')u\,dxdt+
\int_\Gamma\left.\delta^h(w*\rho)u\right|_{t_1}^{t_2}\,dx =
-\int_{t_1}^{t_2}\int_\Omega\langle\nabla\delta^h\widetilde
w,\nabla w\rangle\,d\overline xdt,
\label{eq:primera}\end{equation}
where $\rho'=d\rho/dt$.

The left-hand side mimicks $\iint (\partial_t w)(\partial_t u)dxdt=c\iint
(\partial_t(u^{(m+1)/2}))^2dxdt$, while the right-hand side mimicks $-\iint
\langle\nabla(\partial w/\partial t),\nabla w\rangle\,d\overline xdt $.
Here the times $t_1$ and $t_2$ are subject to be moved
slightly on the condition that $t_1<T_1<T_2<t_2$ for some fixed $0<T_1<T_2$. We
now analyze the different integrals.

$\bullet$  Using \eqref{convol}, the integral in the right-hand side can be
written as
$$
\begin{array}{l}
 -\displaystyle\int_{t_1}^{t_2}\int_\Omega\langle\delta^h\nabla\hat
w,\nabla \hat w\rangle\,d\overline xdt= \\ [4mm]
\displaystyle \frac1{2h}\int_{t_1}^{t_2}\int_\Omega\langle\nabla\hat
w(t-h),\nabla
\hat w(t)\rangle\,d\overline xdt-\frac1{2h}\int_{t_1}^{t_2}\int_\Omega\langle\nabla\hat
w(t+h),\nabla
\hat w(t)\rangle\,d\overline xdt= \\ [4mm]
\displaystyle \frac1{2h}\int_{t_1}^{t_1+h}\int_\Omega\langle\nabla\hat
w(t-h),\nabla
\hat w(t)\rangle\,d\overline xdt-\frac1{2h}\int_{t_2-h}^{t_2}\int_\Omega\langle\nabla\hat
w(t+h),\nabla
\hat w(t)\rangle\,d\overline xdt= \\ [4mm]
\displaystyle \frac1{2h}\int_{t_1-h}^{t_1}\int_\Omega\langle\nabla\hat
w(t+h),\nabla
\hat w(t)\rangle\,d\overline xdt-\frac1{2h}\int_{t_2}^{t_2+h}\int_\Omega\langle\nabla\hat
w(t-h),\nabla
\hat w(t)\rangle\,d\overline xdt.
\end{array}
$$
These two last integrals are bounded, in absolute value, by
$$
\displaystyle\frac1{2h}\int_{t_1-h}^{t_1+h}\int_\Omega|\nabla\hat
w(t)|^2\,d\overline xdt+ \frac1{2h}\int_{t_2-h}^{t_2+h}\int_\Omega|\nabla\hat
w(t)|^2\,d\overline xdt=Y_1
$$
Since $|\nabla w|^2\in L^1$, by picking good times $t_1$ and $t_2$ from a
diadic division of intervals around $T_1$ and $T_2$, the quantity $Y_1$ is
bounded uniformly in $h$ (though $t_1$ and $t_2$ depend on $h$).

$\bullet$ We turn to the left-hand side terms in \eqref{eq:primera}. Assume
first $m>1$. The second integral in \eqref{eq:primera} is bounded using
Proposition~\ref{pro:estimates}.

As to the first integral in  \eqref{eq:primera} we have
$$
\begin{array}{l}\displaystyle-\int_{t_1}^{t_2}\int_\Gamma\delta^h(w*\rho')u\,dxdt=
\int_{t_1+h}^{t_2-h}\int_\Gamma(w*\rho')\delta^hu\,dxdt \\ [3mm]
\qquad\displaystyle+
\frac1{2h}\int_{t_1-h}^{t_1+h}\int_\Gamma(w*\rho')u(t+h)\,dxdt-
\frac1{2h}\int_{t_2-h}^{t_2+h}\int_\Gamma(w*\rho')u(t-h)\,dxdt.
\end{array}
$$
The last two integrals are bounded in absolute value, in the same way as the
second integral in  \eqref{eq:primera}. The first term in the right-hand side
above is the one we want to estimate carefully.  By Lemma~\ref{lem-a2} and
under the extra assumption of monotonicity in time, $\partial_t u\ge
0$, we have
\begin{equation}
I=(w*\rho')(x,t)\ge c\delta^hw(x,t),
\label{wrho}\end{equation}
with a positive constant $c$ that depends on $m,N$ and the smoothing kernel
$\eta_1$. It only remains to observe that
$$
(\delta^hw)\,(\delta^hu)=(\delta^hu^m)\,(\delta^hu)\ge c
(\delta^hu^{(m+1)/2})^2.
$$
This turns out to be an easy calculus  problem, using  the technical Lemma~\ref{lem-a3}.

$\bullet $ Summing up, we get a uniform estimate
$$
\displaystyle\int_{T_1}^{T_2}\int_\Gamma (\delta^h u^{(m+1)/2})^2\, dxdt\le C.
$$
Letting $h\to0$, we get $\partial_t(u^{\frac{m+1}2})\in
L^2(\mathbb{R}^N\times[T_1,T_2])$ for every $0<T_1<T_2$.

$\bullet $  General situation when $m>1$:   we want to apply
Lemma~\ref{lem-a2} to $I=(w*\rho')(x,t)$ without the extra assumption of
monotonicity in time, only using that  if $m>1$, $\partial_t u\ge
-cu$. Then we use the second version of Lemma ~\ref{lem-a2}. An extra term
appears but it is controllable.

$\bullet $  In the case $m<1$ we need to perform some little extra
calculations. First, in order to use a bound of the form \eqref{utmeasure}
for $w$ we take advantage of the fact that $w\ge c>0$ in every compact set of
$\overline\Omega$. Thus we consider a new test function by multiplying
$\varphi$ by a cutoff function $\psi(\overline x)$. The extra terms obtained
are easily bounded.

On the other hand, formula
\eqref{wrho} holds with reverse inequality,
\begin{equation}
I=(w*\rho')(x,t)\le c\delta^hw(x,t),
\label{wrho2}\end{equation}
provided $\partial_t w\le0$, as before. Care has to be taken in the
general case, where we use the estimate $\partial_t w\le cw$.
\qed

We now prove the main result.

\begin{Theorem}
The weak (nonnegative) solution $u$ to Problem \eqref{eq:main} constructed  in
Theorem {\rm\ref{th:main}}  is a strong solution. Moreover,
\begin{equation}
\left\|\dfrac{\partial u}{\partial t}(\cdot,t)\right\|_{L^1(\mathbb{R}^N)}
\le\frac{2}{|m-1|t}\,\|f\|_{L^1(\mathbb{R}^N)}\,.
\label{est-strong}\end{equation}
\end{Theorem}

\noindent{\it Proof. }
Let $w$ be the weak solution to Problem \eqref{pp:local} associated to $u$. We
want to prove that the time derivative of $u$ is actually an integrable
function and that the normal derivative of $w$ on $\Gamma$ is a distribution.
Thus the second equation in \eqref{pp:local}, i.e.~the equation
\eqref{eq:main}, holds almost everywhere.

To deal with the normal derivative  we only have to take into account that,
thanks to the trace embedding, since for every $t>0$ we have $w\in
H^1(\Omega)$, then $\partial_y w(\cdot,0,t)$ is a
distribution in $H^{-1/2}(\Gamma)$.

We now look at $\partial_t u (\cdot,t)$. We use a
technical result by B\'enilan \cite{Benilan2}, see also
\cite{Vazquez:Book}. As in the previous proof, we may assume
$\partial_t u \ge0$. Also we know that $u\in
L^1(\mathbb{R}^N)\cap L^\infty(\mathbb{R}^N)$ for every $t>0$ and
finally, from Proposition \eqref{pro:estimates} we have
$\partial_t  ( u^{(m+1)/2}) \in L^2(\mathbb{R}^N)$.

All these estimates allow us to apply Lemma~8.2 in
\cite{Vazquez:Book} to get $\partial_t u \in
L^p_{loc}(\mathbb{R}^N)$ for every $p\in[1,p_1)$, where
$p_1=\min\{(m+1)/m),2\}$. This gives $\partial_t u \in
L^1_{loc}(\mathbb{R}^N)$ for every $t>0$. Formula
\eqref{est-strong} follows from \eqref{utmeasure}, and in fact $\partial_t u \in
L^1(\mathbb{R}^N)$ for every $t>0$.
 \qed

Once we know that our solution is a strong solution, we can perform the
above-mentioned formal calculation, analogous to the case of the PME.

\begin{Proposition}
In the above hypotheses,
\begin{equation}
\int_{t_1}^{t_2}\int_\Gamma\left|\frac{\partial u^{(m+1)/2}}{\partial
t}\right|^2(x,t)\,dxdt\le
\frac{m+1}{8mt_1}\int_{\Gamma}u^{m+1}(x,t_1)\,dx.
\label{m+1-2}\end{equation}
\label{pro:estim-ut}\end{Proposition}

\medskip

We end this section with two more estimates that will be useful in the sequel.
By comparison assume $f\ge0$. Putting $\varphi=w$ as test function and obtain
\begin{equation}
\int_0^t\int_{\Omega}|\nabla w|^2\,d\overline
xds+
\frac1{m+1}\int_{\Gamma}u^{m+1}\,dx=\frac1{m+1}\int_{\Gamma}f^{m+1}\,dx.
\label{m+1}\end{equation}

We have thus a control of the $L^2$ norm of the gradient in terms of the
initial data.

\begin{Proposition}
  In the above hypotheses
\begin{equation}
\|\nabla w\|_{L^2(\Omega\times(0,\infty))}\le
c\|f\|_{L^{m+1}(\Gamma)}^{(m+1)/2}.
\label{L2grad-m+1}\end{equation}

\end{Proposition}

Another consequence of \eqref{m+1} is that the norm
$\|u(\cdot,t)\|_{L^{m+1}(\mathbb{R}^N)}$ is nonincreasing in time. In fact this
also follows from the elliptic estimates of Section~\ref{sect-exist}.

An easy generalization shows the following property.
\begin{Proposition}
  In the above hypotheses any $L^p$ norm of the solution is nonincreasing in time for every $p\ge
1$. Even more, if $\Psi$ is any convex nonnegative real function, then $
\int_{\mathbb{R}^N}\Psi(u(x,t))\,dx$ is a nonincreasing function.
\label{pro:decr}\end{Proposition}

\noindent{\it Proof. }
We have
$$
\int_\Gamma\Psi(u(x,t_2))\,dx-
\int_\Gamma\Psi(u(x,t_1))\,dx=
-\int_{t_1}^{t_2}\frac 1m\int_{\Omega}w^{(1-m)/m}\Psi''(w^{1/m})|\nabla
w|^2\,d\overline xds\le0.
$$
\qed

\section{Solutions with data in $L^1$}\label{sect-L1}
\setcounter{equation}{0}

In order to construct solutions with initial data in $L^1$ we approximate by
problems with data in $L^1\cap L^\infty$ and use the $L^1$ contractivity to pass
to the limit (together with the estimates of the gradients in $L^2$).

We begin by proving the bound in $L^\infty$ in terms of the $L^1$ norm of the
initial datum, estimate \eqref{eq:L-inf}.

\subsection{Smoothing effect}

\noindent{\it Proof of Theorem \ref{th:L-inf}.}
{\sc Case $m>1$}. We use a technique inspired by the proof in \cite{HP},
Lemma~3.3. To do that we consider the solution $w$ to the local Problem
\eqref{pp:local}. Using estimate~\eqref{ut+u>} and the integral formulation of
that problem, identity \eqref{strong-local}, we get for any $t>0$,
$$
\int_\Omega\langle\nabla w,\nabla\varphi\rangle-L\int_\Gamma u\varphi\le0,
$$
where $L=1/((m-1)t)$, $\varphi\ge0$. The letter $c$ will denote any
constant depending only on $m$ and $N$. Let $\varphi=w^p$, $p>0$. Manipulating
the integral in $\Omega$ and using the trace embedding, we obtain
\begin{equation}
L\int_\Gamma u^{mp+1}\ge\int_\Omega\langle\nabla w,\nabla w^p\rangle=
\frac{4p}{(p+1)^2}\int_\Omega|\nabla w^{\frac{p+1}2}|^2 \ge
c\left(\int_\Gamma u^{\frac{Nm(p+1)}{N-1}}\right)^{\frac{N-1}N}.
\label{uno-1}\end{equation}
Observe that $4p/((p+1)^2)\le1$, so $c$ does not depend on $p$. Therefore, if
we denote $\|u(\cdot,t)\|_{L^q(\Gamma)}$, for the fixed time $t$  by $\|u\|_q$,
we have
$$
\|u\|_{sm(p+1)}\le (cL)^{\frac1{m(p+1)}}\|u\|_{mp+1}^{\frac{mp+1}{m(p+1)}},
$$
where $s=N/(N-1)$. We now iterate this estimate,
$$
\|u\|_{q_{k+1}}\le (cL)^{\frac s{q_{k+1}}}\|u\|_{q_k}^{\frac{sq_k}{q_{k+1}}},
$$
where $q_{k+1}=s(q_k+m-1)$. We obtain
$$
\|u\|_{q_{k}}\le (cL)^{a_k}\|u\|_{q_0}^{b_k},
$$
with the exponents
$$
a_k=\frac1{q_k}\sum_{j=1}^k s^j,\qquad b_k=\frac{s^kq_0}{q_k}\,.
$$
If we start with $q_0>1$, it is a calculus matter to obtain
$$
q_k=As^k-B,\qquad B=\frac{s(m-1)}{s-1}=N(m-1),\quad A=q_0+B,
$$
and then
$$
 \lim_{k\to\infty}b_k=\frac{q_0}A,\quad \lim_{k\to\infty}a_k=\frac NA.
$$
We conclude that
$$
\|u\|_\infty\le cL^{\frac NA}\|u\|_{q_0}^{\frac{q_0}A}.
$$
Now recall that $q_0>1$. In order to reach the  $L^1$ norm, we use a classical
interpolation argument once the above partial smoothing effect is obtained. To
do this we put $q_0=2$ for simplicity, since any other value gives the same
estimate. To prove an $L^1$-$L^2$ smoothing effect we observe that
$$
\|u\|_2\le \|u\|_1^{1/2}\|u\|_\infty^{1/2}\le
cL^{\frac N{2mA}}\|u\|_1^{\frac12}\|u\|_2^{\frac1{mA}},
$$
where here $A=2+N(m-1)$. Therefore
$$
\|u\|_2\le
L^{\frac{N}{2(1+N(m-1))}}\|u\|_1^{\frac{2+N(m-1)}{2(1+N(m-1))}},
$$
which yields the estimate
$$
\|u\|_\infty\le c
L^{\frac{N}{1+N(m-1)}}\|u\|_1^{\frac1{1+N(m-1)}}.
$$
Recalling that $L\sim 1/t$, and using the Remark after Theorem~\ref{th:plocal}, this gives
\eqref{eq:L-inf}.

{\sc Case $m\le 1$}. Since we do not have a lower estimate for $\partial_tu$ in
this case, we use here a different idea: we look at the proof of the smoothing
effect for the quasi-geostrophic equation performed in
\cite{CV}, see also \cite{CW}. In fact, the proof can be adapted to get also an
$L^{1+\varepsilon}$-$L^{\infty}$ regularizing effect, for any $\varepsilon>0$.
In order to reach $L^1$ we use an analogous procedure as before, but here the
interpolation must be performed iteratively, since the first regularizing
effect in this case is obtained for different times.

Let $w$ be the solution to Problem \eqref{pp:local}, and fix $C>0$ a constant
to be chosen later. Consider the functions $
\zeta=(w-C)_+$ and $\theta=(w^{1/m}-C^{1/m})_+$. By Kato's inequality we have that $\zeta$ is subharmonic.
Using  $\varphi(\theta)$ as test function in the definition of Problem
\eqref{pp:local},
we get
$$
\frac d{dt}\int_\Gamma
\psi(\theta)\,dx+\int_\Omega \langle\nabla\varphi(\theta),\nabla\zeta\rangle\,d\overline
x\le0,
$$
where $\psi'=\varphi$. If we put $\varphi(\theta)=\theta^\varepsilon$, the
following will produce an $L^{1+\varepsilon}$-$L^{\infty}$ regularizing effect.
As we have said before, it is enough to get an $L^2$-$L^\infty$ regularizing
effect in the first step. Thus we consider for simplicity
$\varphi(\theta)=\theta$. We have then
$$
\frac d{dt}\int_\Gamma
\theta^2\,dx+2\int_\Omega \langle\nabla\theta,\nabla\zeta\rangle\,d\overline
x\le0
$$
We now fix the constants. Let $t_0>0$, $t_k=t_0(1-2^{-k})$, $M>0$,
$C=C_k=(M(1-2^{-k}))^m$, $\theta=\theta_k$, $\zeta=\zeta_k$. Define the
quantity,
\begin{equation}
U_k=\sup_{t\ge t_k}\left( \int_\Gamma\theta_k^{2}(x,0,t)dx\right)+2
\int_{t_k}^\infty\int_{\Omega}\langle\nabla\theta_k,\nabla\zeta_k\rangle\,d\overline
xdt\equiv I_k+Y_k. \label{uk}\end{equation} Our purpose is to obtain a
recursive estimate of $U_k$ in order to prove $\lim\limits_{k\to\infty}U_k=0$.
This will mean $\lim\limits_{k\to0}(\Tr(w^{1/m})-C_k^{1/m})_+=0$, i.e. $u\le M$.
We must choose $M$ carefully.

First we deduce the upper estimate, for $t_{k-1}<s<t_k$,
$$
U_k\le2\int_\Gamma\theta_k^{2}(x,0,s)dx.
$$
Taking the mean over the interval $[t_{k-1},t_k]$ we get
\begin{equation}
U_k\le\frac{2^{k+1}}{t_0}\int_{t_{k-1}}^\infty\int_\Gamma\theta_k^{2}(x,0,s)dxds.
\label{uk-upper}\end{equation} See \cite{CV}. To get a lower
estimate, in the same way as in \cite{CV}, we have to take care of the
exponents. Observe that the spatial integral in $Y_k$ can be written as
$$
\frac{8m}{(m+1)^2}
\int_\Omega|\nabla(w^{\frac{m+1}{2m}}-C_k^{\frac{m+1}{2m}})_+|^2\,d\overline
x.
$$
Using the trace embedding, and the fact that $m<1$, we get
$$
\begin{array}{l}
\displaystyle\int_\Omega|\nabla(w^{\frac{m+1}{2m}}-C_k^{\frac{m+1}{2m}})_+|^2\,d\overline
x\ge
c\left(\int_\Gamma|(w^{\frac{m+1}{2m}}-C_k^{\frac{m+1}{2m}})_+|^{\frac{2N}{N-1}}\,d
x\right)^{\frac{N-1}N}\\ [5mm] \ge
c\displaystyle\left(\int_\Gamma|(w^{\frac{1}{m}}-C_k^{\frac{1}{m}})_+|^{\frac{N(m+1)}{N-1}}\,d
x\right)^{\frac{N-1}N} =
c\left(\int_\Gamma\theta_k^{\frac{N(m+1)}{N-1}}\,dx\right)^{\frac{N-1}N},
\end{array}$$
where $c$ depends only on $m$ and $N$. From this estimate we obtain, using
Riesz interpolation,
$$
\begin{array}{rl}
U_k&\ge\displaystyle I_k^\nu Y_k^{1-\nu}\ge c\sup_{t\ge t_k}\left(
\int_\Gamma\theta_k^{2}\right)^\nu\,\left(\int_{t_k}^\infty\left(
\int_\Gamma\theta_k^{\frac{N(m+1)}{N-1}}\right)^{\frac{N-1}N}\right)^{1-\nu}
\\ [5mm] &\ge
\displaystyle c\left[
\int_{t_k}^\infty\left(\int_\Gamma\theta_k^{2}\right)^{\frac{\nu}{1-\nu}}\,
\left(\int_\Gamma\theta_k^{\frac{N(m+1)}{N-1}}\right)^{\frac{N-1}N}\right]^{1-\nu}\ge
\displaystyle c\left[
\int_{t_k}^\infty\int_\Gamma\theta_k^q\right]^{1-\nu},
\end{array}$$
by taking $\nu=1/(N+1)$, $q=m+1+2/N$. Observe now that $q>2$
provided $m>1-2/N$, which is guaranteed by the condition $m>m_*$. The last step
is to use the above inequality to deduce the required iteration
$$
\begin{array}{rl}
U_k&\le
\displaystyle\dfrac{2^{k+1}}{t_0}\int_{t_{k-1}}^\infty\int_\Gamma\theta_k^{2}(x,0,s)dxds
\\ &\le
\displaystyle\dfrac{2^{k+1+k(q-2)}}{t_0M^{q-2}}
\int_{t_{k-1}}^\infty\int_\Gamma\theta_{k-1}^{q}(x,0,s)dxds \le
\frac{c2^{k(q-1)}}{t_0M^{q-2}} U_{k-1}^{1+1/N}.
\end{array}$$
This gives the limit mentioned above, $\lim\limits_{k\to\infty}U_k=0$, provided
$$
\frac{cU_0^{1/N}}{t_0M^{q-2}}<1.
$$
This can be achieved choosing $M=c t_0^{-\mu}U_0^{\mu/ N}$ large enough, where
$1/\mu=q-2=m-1+2/N$. This gives $u\le M$. Since $U_0\le c\|f\|_{2}^{2}$, we
obtain the $L^2$-$L^\infty$ estimate
\begin{equation}
|u(x,t)|=\le ct^{-\mu}\|f\|_{2}^{2\mu/N}.
\label{firstmu}\end{equation}

We now perform the interpolation. Using \eqref{firstmu} in the interval
$[t/2,t]$, and then in the interval $[t/4,t/2]$, we get
$$
\begin{array}{rl}
\|u(\cdot,t)\|_\infty&\le c(t/2)^{-\mu}\|u(\cdot,t/2)\|_2^{2\mu/N}\le c2^\mu
t^{-\mu}\|u(\cdot,t/2)\|_1^{\mu/N}\|u(\cdot,t/2)\|_\infty^{\mu/N} \\ [4mm]&\le
c2^\mu
t^{-\mu}\|u(\cdot,t/2)\|_1^{\mu/N}\Big(c(t/4)^{-\mu}\|u(\cdot,t/4)\|_2^{2\mu/N}\Big)^{\mu/N}.
\end{array}$$
Iterating this calculation in intervals of the form $[2^{-k}t,2^{-(k-1)}t]$,
using Proposition~\ref{pro:decr}, we obtain
$$
\|u(\cdot,t)\|_\infty\le c^{a_k}2^{b_k}
t^{-d_k}\|u(\cdot,0)\|_1^{e_k}\|u(\cdot,2^{-k}t)\|_2^{f_k}.
$$
Using the fact that  $m>m_*$ implies $\mu<N$, we see that the exponents
satisfy, in the limit $k\to\infty$,
$$
\begin{array}{l}
\displaystyle
a_k=\sum_{j=0}^{k-1}\Big(\frac\mu
N\Big)^j\to\frac{(m-1)N+2}{(m-1)N+1}=\frac\gamma N+1\,, \\ [4mm]
\displaystyle
b_k=\sum_{j=0}^{k-1}\mu (j+1)\Big(\frac\mu N\Big)^j\to b<\infty\,, \\ [4mm]
\displaystyle
d_k=\mu a_k\to\frac{N}{(m-1)N+1}=\gamma\,,
\\ [4mm]
\displaystyle
e_k=a_k-1\to \frac\gamma N\,,
\\ [4mm]
\displaystyle
f_k=2\Big(\frac\mu N\Big)^k\to 0\,.
\end{array}
$$
This ends the proof. \qed

\noindent{\it Remark. } Once we assert that the constant $c$ in \eqref{eq:L-inf} is
universal, then the values of the exponents $\gamma$ and $\gamma/N$ are given
as an immediate consequence of the invariance of the equation under the
two-parameter scaling group. This is similar to what happens in the Sobolev
inequalities, cf.~\cite{Evans}, page 262, or what happens in the  PME,
cf.~\cite{JLVSmoothing}, page 29.

\subsection{$L^1$ contraction}

\noindent{\it Proof of Theorem \ref{th:contract}.}
We first recall that $L^1$-energy solutions are weak energy solutions for every
$\tau>0$. Therefore, by the trace embedding we have $u(\cdot,\tau)\in
L^p(\mathbb{R}^N)$ for $1\le p\le 2N/(N-1)$ if $N>1$ (for every $p\ge1$ if
$N=1$). Now condition $m>m_*$ implies $u(\cdot,\tau)\in L^{m+1}(\mathbb{R}^N)$.
This allows us to use the proof of uniqueness for weak energy solutions, cf.
Lemma~\ref{oleinik} and the Remark after it, for every $\tau>0$. Continuity in
$L^1$ gives uniqueness up to $\tau=0$. We now approximate the initial datum
$f\in L^1(\mathbb{R}^N)$ by functions $f_k\in L^1(\mathbb{R}^N)\cap
L^\infty(\mathbb{R}^N)$, $f_k\to f$ in $L^1(\mathbb{R}^N)$. The above implies
that corresponding solutions $u_k$ converge in $L^1(\mathbb{R}^N)$ to our
solution $u$, and the contractivity result in \eqref{th:plocal} passes to the
limit.
\qed

\subsection{End of proof of Theorem \ref{th:main2}}

Once we have the $L^1$ contraction \eqref{eq:contract-l1} at our disposal, and
using the above approximation of $u$ by $u_k$, we observe that the estimate
\eqref{eq:L-inf} does not depend on the $L^\infty$ norm of the approximations,
and thus it is true for the limit.

As to the conservation of mass, we recall that the proof of
Theorem~\ref{th:mass} relies on bounds for the $L^p$ norms of $u(\cdot,t)$ for
$t>0$. This is handled with the smoothing effect just proved. We thus have,
repeating the proof of Theorem~\ref{th:mass},
$$
\left|\int_{\mathbb{R}^N}\Big(u(x,t)-f(x)\Big)\varphi_R(x)\,dx\right|\le
cR^{-2}\int_0^t\int_{R<|\overline x|<2R}|w(\overline x,s)|\,d\overline xds.
$$
We now use \eqref{eq:L-inf} to estimate the right-hand side $I$. If $m>1$ we
get
$$
I\le
cR^{-1}\int_0^t\|u(\cdot,s)\|^{m-1}_{L^\infty(\mathbb{R}^N)}\|u(\cdot,s)\|_{L^1(\mathbb{R}^N)}\,ds
\le cR^{-1}\|f\|_1^{\frac{\gamma(m-1)}N+1}\int_0^ts^{-\gamma(m-1)}\,ds\to0,
$$
since $\gamma(m-1)<1$. In the case $m<1$ we obtain, instead,
$$
\begin{array}{rl}
\displaystyle I &\le
\displaystyle cR^{N(p-1)/p-1}\|u(\cdot,t)\|_{L^\infty(\mathbb{R}^N)}^{m-1/p}\|u(\cdot,t)\|_{L^1(\mathbb{R}^N)}^{1/p}
\\ [3mm]
&\displaystyle \le
cR^{N(p-1)/p-1}\|f\|_1^{\frac{\gamma(m-1/p)}N+\frac1p}\int_0^ts^{-\gamma(m-1/p)}\,ds\to0,
\end{array}
$$
by choosing $1/m<p<N/(N-1)$, which is possible whenever $m>m_*$.

Finally, the proof of positivity and regularity in Theorems~\ref{th:positivity}
and~\ref{th:regularity2} also apply using again the smoothing effect.
\qed

\section{Continuous dependence}\label{sect.cont.dep.}
\setcounter{equation}{0}

The aim of this section is to prove the continuous dependence of the solutions
constructed in this paper with respect to the initial data $f$ and exponent
$m$. This is true for $m>m_*$.  Let us introduce the notation: if $u=u_{m,f}$ is the
solution corresponding to $m$ and $f$, we write $S(m,f)=u_{m,f}$.
\begin{Theorem}
The map $S: (m_*,\infty)\times L^1(\mathbb{R}^N)\to
C([0,T]:\,L^1(\mathbb{R}^N))$ is continuous. \label{th:cont}
\end{Theorem}

This  will follow from a
result of nonlinear Semigroup Theory which states that if each of $A_n$,
$n=1,2,\dots,\infty$ is an $m$-accretive operator in a Banach space $\cal X$,
$f_n\in\overline{D(A_n)}$ and $u_n$ is the solution of
$$
\frac{du_n}{dt}+A_nu_n=0,\qquad u_n(0)=f_n,
$$
then $A_n\to A_\infty$, $f_n\to f_\infty$ implies $u_n\to u_\infty$ in
$C([0,\infty):\, {\cal X})$, where $A_n\to A$ means
$$
\lim_{n\to\infty}(I+A_n)^{-1}g= (I+A_\infty)^{-1}g\quad\mbox{for all
} g\in {\cal X}.
$$
See, e.g, \cite{Crandall}, \cite{Evans} for statements and references. Hence, the theorem will be a corollary of the convergence of $(I+A_{m_n})^{-1}$, where $A_{m}(u)=(-\Delta)^{1/2}u^{m}=-\partial_y w$. This is what we prove next.

\begin{Proposition}
Let $\{m_n\}_{n=1}^\infty$ be a bounded sequence of numbers in $(m_*,\infty)$
such that $\lim\limits_{n\to\infty} m_n=\bar m>m_*$. Then
$\lim\limits_{n\to\infty}(I+A_{m_n})^{-1}g= (I+A_{\bar m})^{-1}g$ for all $g\in
L^1(\mathbb{R}^N)\cap L^\infty(\mathbb{R}^N)$.
\end{Proposition}
\noindent{\it Proof. } We borrow ideas from \cite{BEG}. Note that an analogous result for the PME  was proved in \cite{Benilan-Crandall}. Let $u_{m}=(I+A_{m})^{-1}g$. The $L^1$-contraction
estimate~\eqref{eq:contract-e} implies the bounds
\begin{equation}
\label{eq:resolvents}
\begin{array}{l}
\|u_{m_n}\|_{L^1(\mathbb{R}^N)}\le \|g\|_{L^1(\mathbb{R}^N)}\,,\\
[3mm] \|u_{m_n}-\tau_h
u_{m_n}\|_{L^1(\mathbb{R}^N)}\le\|g-\tau_hg\|_{L^1(\mathbb{R}^N)}\,,
\end{array}
\end{equation}
for each $h\in\mathbb{R}^N$, where $(\tau_h v)(x)=v(x+h)$. This is
enough, thanks to Fr\'echet-Kolmogorov's compactness criterium
\cite{Brezis2}, to prove that $\{u_{m_n}\}$ is precompact in
$L^1(K)$ for each compact set $K\subset\mathbb{R}^N$.

To extend compactness to the whole $\mathbb{R}^N$ we need to control
the tails of the solutions uniformly. More precisely, we need to
prove that, given $\epsilon>0$, there exists some $R>0$ such that
$\|u_{m_n}\|_{L^1(\mathbb{R}^N\setminus B_{R}(0))}<\epsilon$. This
follows from a computation which is very similar to the one in the
proof of Theorem~\ref{th:mass}, but now taking as test function
$1-\varphi_R$ instead of $\varphi_R$. First observe that
$u_{m_n}=(I+A_{m_n})^{-1}g$ means
$$
\int_\Gamma u_{m_n}\varphi_R\,dx=\int_\Gamma
g\varphi_R\,dx-\int_\Omega\langle\nabla
w_{m_n},\nabla\varphi_R\rangle\,d\overline x.
$$
Therefore, with the above mentioned test function we obtain
$$
\int_{\{|x|>R\}}u_{m_n}\,dx\le\int_{\{|x|>R\}}g\,dx+cR^{-2}\int_{\{R<|\overline x|<2R\}}w_{m_n}\,d\overline
x\rightarrow0
$$
as $R\to\infty$ uniformly in $n$, see the proof of Theorem~\ref{th:mass}.

We have obtained that along some subsequence, which we also call $\{m_n\}$, the
following convergence holds
$$
u_{m_n}\rightarrow\;u_*,\qquad \mbox{in } L^1(\Gamma),
$$
for some function $u_*$. Using the Poisson kernel, we also have
$w_{m_n}\to w_*$ in $L^1(\Omega)$, where $w_*$ is the harmonic
extension of $u_*$ to the upper half-space. On the other hand, we
have a uniform control in  $L^2(\Omega)$ of the gradients of $w$ in
terms of the data $g$, see~\eqref{rabia}. Thus, there is weak
convergence in $L^2(\Omega)$ of the gradients $\nabla w_{m_n}$ along
some subsequence towards  $\nabla w_*$. All this is enough to pass
to the limit in~\eqref{eq:resolvents} to show that the limit $u_*$
is indeed $u_{\overline m}$. \qed

\noindent{\it Remark. } If $m>1$, an easier alternative proof can be
performed using the compactness results of \cite{Simon}.

\section{Comments and extensions}\label{sec.cae}

\noindent {\sc Alternative approaches.}~It is not difficult to prove a posteriori that
the constructed semigroup is also contractive with respect to the norm
$H^{-1/2}$. This property could be used as a starting point in the existence
and uniqueness theory by using results of the theory of monotone operators in
Hilbert spaces, as developed in \cite{Brezis} for the PME case. We have chosen
our present formulation because we have found a number of advantages in proceeding in this manner.

On the other hand, Crandall and Pierre  developed in \cite{CP} an
abstract approach to study evolution equations of the form
$\partial_t u+A\varphi(u)=0$ when $A$ is an $m$-accretive operator
in $L^1$ and $\varphi$ is a monotone increasing real function.  This allows
to obtain a  \emph{mild} solution using the Crandall-Liggett Theorem. Our problem falls within
this framework. Let us point out that such an abstract construction does not give
enough information to prove that the mild solution is in fact a weak
solution, in other words, to identify the solutions in a differential sense.

\medskip

\noindent {\sc Extension.}~As a natural extension of this work we can consider the more general
model based on the equation
\begin{equation}\frac{\partial u}{\partial t} +
(-\Delta)^{\sigma/2} (|u|^{m-1}u)=0,
\end{equation}
where the fractional Laplacian has an exponent $\sigma\in (0,2)$, and $m>m_*$
for some $m_*(N,\sigma)\in [0,1)$. Though the main qualitative results are
similar to the ones presented here, the theory of these fractional
operators with $\sigma\ne 1$ has some technical difficulties that make it
convenient to be treated at a second stage. We recall that Caffarelli and Silvestre
\cite{CS} have recently characterized the \emph{Laplacian of order} $\sigma$,
$(-\Delta)^{\sigma/2}$, by means of another auxiliary extension approach. We will
use such an extension in a separate paper, \cite{PQRV}, to treat in detail the more general
fractional diffusion model.

\section{Appendix}\label{appendix}
\setcounter{equation}{0}

We recall some technical results that we have needed in the proof of the property of strong solutions.
The first seems to be well-known.

\begin{Lemma} On the condition that $\eta$ is a smooth convolution kernel with $\eta(-x)=\eta(x)$ we have
 for every pair of $L^2$ functions in $\mathbb{R}$
 \begin{equation}
\int (h*(\eta*\eta))(t)g(t)\,dt=\int (h*\eta)(t)(g*\eta)(t)\,dt.
\label{convol-a1}\end{equation}

\label{lem-a1}\end{Lemma}

\begin{Lemma} $i)$ Let $g$ be a positive nondecreasing function. If $\rho_h(t)=(1/h)\rho_1(t/h)$, where $\rho_1$
is a smooth, symmetric, compactly supported nonnegative function, with support
$[-2,2]$, and $\delta^hg(t)=(g(t+h)-g(t-h))/(2h)$, then
  \begin{equation}
(g*\rho_h')(t)\ge c\delta^hg(t).
\label{convol-a2}\end{equation}
\noindent $ii)$ If instead of nondecreasing we have  $g'(t)\ge -Ag(t)$, the conclusion is
 \begin{equation}
(g*\rho_h')(t)\ge c\delta^hg(t)- cAg(t).
\label{convol-a2+}\end{equation}
\label{lem-a2}\end{Lemma}

\noindent{\it Proof. } $i)$ In the case of nondecrasing $g$ we do as follows:
$$
\begin{array}{rl}
(g*\rho_h')(t)&\displaystyle=\int_0^{2h}\Big(g(t+s)-g(t-s)\Big)(-\rho'(s))\,ds \\
[3mm] &\displaystyle\ge\int_h^{2h}\Big(g(t+h)-g(t-h)\Big)(-\rho'(s))\,ds \\
[3mm]
&\displaystyle\ge\delta^hg(t)\int_h^{2h}(-h\rho'(s))\,ds=\delta^hg(t)\rho_1(1).
\end{array}
$$
We pass from the fist to the second line using the positivity of the integrand and the fact that for
$s\in (h,2h)$ we have $g(t+s)\ge g(t+h)\ge g(t-h)\ge g(t-s)$.

\medskip

\noindent $ii)$ The difference is now that we have to use the weaker inequality $g(t+s)\ge g(t)e^{-As}$ if $s>0$.
We have again
$$
(g*\rho_h')(t) \displaystyle=\int_0^{2h}\Big(g(t+s)-g(t-s)\Big)(-\rho'(s))\,ds.
$$
Since $g(t+s)-g(t-s)\ge -g(t)(e^{As}-e^{-As})$ we get
$$
\int_0^{h}\Big(g(t+s)-g(t-s)\Big)(-\rho'(s))\,ds\ge -g(t)\int_0^{h}(e^{As}-e^{-As})(-\rho'(s))\,ds
$$
and the last integral is bounded uniformly for $h$ small in the form $cA$. For
the other part of the integral we have
$$
\begin{array}{rl}
&\displaystyle\int_h^{2h}\Big(g(t+s)-g(t-s)\Big)(-\rho'(s))\,ds \\
&\displaystyle\ge\int_h^{2h}\Big(g(t+h)e^{-A(s-h)}-g(t-h)e^{A(s-h)}\Big)(-\rho'(s))\,ds \\
&\displaystyle\ge\delta^hg(t)\int_h^{2h}(-h\rho'(s))\,ds-I_1
= \delta^hg(t)\rho_1(1)-I_1\,,
\end{array}
$$
where
$I_1=\int_h^{2h}\Big(g(t+h)(1-e^{-A(s-h)})+g(t-h)(e^{A(s-h)}-1)\Big)(-\rho'(s))\,ds$
can also be estimated as $cAg(t)$. \qed

We end this list of results with an easy but useful calculus lemma.

\begin{Lemma} There exists a positive constant $c$ depending on $m>0$ such that
  \begin{equation}
(x^m-1)(x-1)\ge c(x^{\frac{m+1}2}-1)^2,\qquad\forall\;x\ge1.
\label{convol-a3}\end{equation}
\label{lem-a3}
and also
 \begin{equation}
(x^m+1)(x+1)\ge c(x^{\frac{m+1}2}+1)^2,\qquad\forall\;x\ge1.
\label{convol-a3+}\end{equation}
\label{lem-a3+}
\end{Lemma}

\noindent{\it Proof. }
The quotient of the two positive  functions $F(x)=(x^m-1)(x-1)$ and
$G(x)=(x^{\frac{m+1}2}-1)^2$ is bounded below away from zero in the interval
$[1,\infty)$ unless the limit at $x=1$ or $x\to \infty$ is zero. It is clear
that at infinity the limit is 1, whereas at $x=1$ we can use L'Hopital to get
$$
\lim_{x\to 1}\frac{F(x)}{G(x)}=\frac{4m}{(m+1)^2}.
$$
This number is positive and equal or less than 1. The other inequality is similar. \qed


\vskip .5cm

\noindent \textsc{Acknowledgment.} All the authors supported by
Spanish Projects MTM2008-06326-C02-01 and -02 and by ESF Programme
``Global and geometric aspects of nonlinear partial differential
equations".  We thank I. Athanasopoulos and L. Caffarelli for
comments on their work \cite{AC}, and M. Pierre for bringing to our
attention  the paper \cite{CP}.



\

\noindent{\bf Addresses:}

\noindent{\sc A. de Pablo: } Departamento de Matem{\'a}ticas, Universidad Carlos III de Madrid, 28911 Legan{\'e}s,
Spain. (e-mail: arturo.depablo@uc3m.es).

\noindent{\sc F. Quir\'{o}s: }
Departamento de Matem\'{a}ticas, Universidad Aut\'{o}noma de Madrid, 28049 Madrid, Spain.
(e-mail: fernando.quiros@uam.es).

\noindent{\sc A. Rodr{\'\i}guez: }
Departamento de Matem\'{a}tica, ETS Arquitectura, Universidad Polit\'{e}cnica de Madrid, 28040 Madrid, Spain. (e-mail: ana.rodriguez@upm.es).

\noindent{\sc J.~L. V{\'a}zquez: }
Departamento de Matem\'{a}ticas, Universidad Aut\'{o}noma de Madrid, 28049
Madrid, Spain. (e-mail: juanluis.vazquez@uam.es). Second affiliation: Institute ICMAT.

\end{document}